\documentclass{amsart} \usepackage{amssymb, amscd}
\numberwithin{equation}{section}

\def\BB{{\mathbb B}}

\def\CC{{\mathbb C}} 
\def\DD{{\mathbb D}}

\def\EE{{\mathbb E}}

\def\LL{{\mathbb L}} 

\def\PP{{\mathbb P}} 
\def\QQ{{\mathbb Q}} 
\def\RR{{\mathbb R}}
\def\ZZ{{\mathbb Z}}

\def\G{{\Gamma}} 
\def\g{{\gamma}}

\def\st{{\rm st}} 
\def\ss{{\rm ss}} 
 
\def\tor{{tor}} 
\def\bb{{bb}}

\def\bs{\backslash} 
\def\bss{{\bs\!\bs}} 
\def\pt{{\bullet}}
\def\one{{\mathbf 1}}

\def\Ccal{{\mathcal C}}

\def\Hcal{{\mathcal H}} 
\def\Ical{{\mathcal I}}

\def\Lcal{{\mathcal L}}

\def\Ocal{{\mathcal O}}

\def\la{\langle} 
\def\ra{\rangle}

\newcommand\Hom{\operatorname{Hom}}

\newcommand\im{\operatorname{Im}}
 
\newcommand\lie{\operatorname{Lie}}
\newcommand\pic{\operatorname{Pic}}

\newcommand\re{\operatorname{Re}}
 
\newcommand\sym{\operatorname{Sym}} 
\newcommand\GL{\operatorname{GL}}
\newcommand\PO{\operatorname{PO}}

\newcommand\SL{\operatorname{SL}} 
\newcommand\SO{\operatorname{SO}}
\newcommand\so{\operatorname{\mathfrak{so}}} 
\newcommand\Orth {\operatorname{O}} 

\newcommand\PSL{\operatorname{PSL}}

\newtheorem{theorem}{Theorem}[section]
\newtheorem{lemma}[theorem]{Lemma}
\newtheorem{lemmadef}[theorem]{Lemma-Definition}
\newtheorem{proposition}[theorem]{Proposition}
\newtheorem{corollary}[theorem]{Corollary}

\theoremstyle{definition} 
\newtheorem{definition}[theorem]{Definition}

\newtheorem{example}[theorem]{Example}
\newtheorem{setting}[theorem]{The setting} 
\newtheorem{none}[theorem]{}

\theoremstyle{remark} 
\newtheorem{remark}[theorem]{Remark}

\newtheorem{problem}[theorem]{Problem}
\newtheorem{notation}[theorem]{Notation}

\begin{document}

\centerline{\bf {COMPACTICATIONS DEFINED BY ARRANGEMENTS II:}}
\title{Locally symmetric varieties of type IV}
\author{Eduard Looijenga}
\thanks{Part of this work was carried out at the Mathematical Sciences Research Institute
at Berkeley  and supported by NSF grant DMS-9810361}

\address{Faculteit Wiskunde en Informatica, Universiteit Utrecht,
Postbus 80.010, NL-3508 TA Utrecht, Nederland}
\email{looijeng@math.uu.nl}



\begin{abstract} 
We define a new class of completions of locally symmetric varieties of type IV
which interpolates between the Baily-Borel compactification and Mumford's toric 
compactifications. An arithmetic arrangement in a locally symmetric variety of 
type IV determines such a completion canonically. This completion admits a natural 
contraction that leaves the complement of the arrangement untouched. The 
resulting completion of the arrangement complement is very much like a Baily-Borel 
compactification: it is the proj of an algebra of meromorphic automorphic forms. 
When that complement has a moduli space interpretation, then
what we get is often a compactification obtained by means of geometric invariant 
theory. We illustrate this with several examples: moduli spaces of polarized $K3$ 
and Enriques surfaces and the semi-universal deformation of a triangle singularity.

We also discuss the question when a type IV arrangement is definable by an 
automorphic form.
\end{abstract}

\maketitle
\section*{Introduction}
There are two classes of irreducible bounded symmetric domains which
admit complex totally geodesic hypersurfaces: the balls, which are associated 
to unitary groups of signature $(1,n)$, and the type IV-domains, 
associated to real orthogonal groups of signature $(2,n)$.  
A domain of either class comes naturally with an  embedding in a complex projective 
space for which its totally geodesic hypersurfaces are precisely the 
hyperplane sections. This is why we call a locally finite collection of 
these an \emph{arrangement} on such a domain. We are interested in the
case when this arrangement is \emph{arithmetic} in the sense that it is
a finite union of orbits of an arithmetic group of automorphisms of the domain 
in the set of hyperplane sections.
If we denote the domain by $\DD$ and the arithmetic group by $\G$, then
$X:=\G\bs\DD$ is a quasi-projective orbifold and our assumption
implies that the arrangement defines a hypersurface $D$ in $X$;  
a Cartier divisor on $X$ (when viewed as an orbifold) supported by $D$
amounts to assigning to every $\G$-orbit of the arrangement a
nonzero integer. We investigated this situation already for ball quotients in 
Part I \cite{I}; here we concentrate on the more involved type IV quotients 
and so in the remainder of this introduction we shall assume that we are in that case.

The variety $X$ has as its natural projective completion the Baily-Borel 
compactification $X^\bb$. If $\dim X\ge 3$, then $X^\bb$ has the short (though 
not very informative) characterization as the proj of the algebra of 
$\G$-automorphic forms on the domain. Its boundary is of dimension $\le 1$
and $X^\bb$ is in general highly singular there---for this reason we cannot
expect the natural extension $D^\bb$ of $D$ to $X^\bb$
to be the support of a Cartier divisor. On the other hand, an 
automorphic form always defines a Cartier divisor on $X^\bb$, so this is 
presumably the main reason why $D$ is usually not definable by an automorphic 
form. Indeed, the necessary and sufficient condition that we find for a divisor 
supported by $D$ to be Cartier shows that this is a highly nontrivial property.
This is also the point of view taken by Bruinier and Freitag in
\cite{bruinier}, who show that in some special cases the Cartier property 
along one-dimensional boundary components suffices for $D$ being
the zero set of an automorphic form that admits a product expansion.
 
But our main goal is different, namely to define a natural projective compactification 
of $X^\circ:=X-D$ in the spirit of the Satake, Baily-Borel theory. 
This compactification, which we denote by $\widehat{X^\circ}$, is obtained
via two intermediate compactifications which have an interest in their own right:
First we define a normal  blowup of $X^\bb$ which leaves $X$ untouched and has
the virtue that the closure of each irreducible component
of $D$ in this blowup is Cartier. This is basically the minimal normal blowup
with this property. We do this in a constructive, 
Satake-like, fashion. The strict transforms of the irreducible components of $D$
form what we still might call an arrangement on this blowup. We blow up 
this arrangement in a rather obvious manner and then we prove that this second
blowup has a natural blowdown---this is our $\widehat{X^\circ}$. The
construction also yields an ample line bundle $\widehat\Lcal$ on 
$\widehat{X^\circ}$ which has the same restriction to $X^\circ$ as the automorphic 
line bundle on  $X$. So every meromorphic $\G$-automorphic form
whose polar locus is contained in the arrangement defines a meromorphic section
of $\widehat\Lcal$ that is regular on $X^\circ$. This  has an interesting consequence, 
as we shall now explain. If $D$ is the zero set of an automorphic form, then 
the boundary of this compactification is  always a hypersurface. But that
situation is rather special  and quite often this boundary has codimension 
$\ge 2$. This implies
that a meromorphic $\G$-automorphic form as above defines a regular section
of a nonnegative tensor power of $\widehat\Lcal$. So the algebra of such forms
is in fact finitely generated with positive degree generators and
$\widehat{X^\circ}$ can be characterized  as the proj of this algebra.
This fact is very useful in the situation where 
the automorphic bundle over $X^\circ$ also appears as the quotient 
of an ample line bundle over a normal variety relative a proper action of
a reductive algebraic group: we show that $\widehat{X^\circ}$ is then likely
to appear as a GIT compactification. This happens for instance for 
the moduli spaces of sextic curves, of quartic surfaces (in both cases
simple singularities allowed) and of curves on $\PP^1\times\PP^1$ 
of bidegree $(4,4)$, which have in 
common that they all parametrize certain families of $K3$-surfaces.
In these cases, the first intermediate 
blowup of $X^\bb$ has also algebro-geometric relevance.
An application in the same spirit concerns the semi-universal deformation of 
a triangle singularity. These applications are all discussed in this paper, 
although we intend to discuss some of them in more detail elsewhere.

Let us now say something about the way this paper is organized. The first two 
sections are essentially a review of the Baily-Borel compactification of a 
locally symmetric variety of
type IV: Section $1$ is devoted to the domains of that type, its boundary
components and relevant groups, whereas Section $2$ deals with the
Baily-Borel compactification  proper. We  use the occasion not only to set up 
notation, but also to state the basic results in a manner that is best 
suited for our purposes. This leads to a purely geometric approach, which
almost avoids any reference to the general theory of algebraic groups 
(for instance, parabolic subgroups are never mentioned).
The next four sections make up what is perhaps the core of this paper: the main goal
here is to define a modification of the Baily-Borel compactification $X^\bb$ in
terms of combinatorial data. We first do this over the one-dimensional boundary
components. Here the situation is relatively simple: the  data in question then
consist of giving a $\QQ$-vector subspace of a certain vector space. 
Such a subspace comes up naturally when we are  given an arithmetic arrangement.
We also give a necessary and sufficient condition for an arrangement 
to support a $\QQ$-Cartier divisor near a one-dimensional boundary stratum. 
At this point it takes us only half a page to prove that the `rational double point locus' 
in the moduli space of semipolarized $K3$-surfaces of genus $g>2$
does not extend as a $\QQ$-Cartier divisor in its Baily-Borel compactification.
This implies that this locus is not definable by an automorphic form, 
thus completing an earlier result due to Nikulin (which says that this is so for infinitely 
many values of $g$).   The same program for the  zero dimensional boundary
strata (the `cusps') is carried out in  Sections $4$ and $5$: the former treats the 
modification of the local Baily-Borel compactification, the latter connects this
with arrangements.  The necessary and sufficient condition that we find for an 
arrangement to support a Cartier divisor near  a cusp is familiar in
the theory of generalized root systems \cite{looij:genroot} and Lie-algebra's
of Kac-Moody type. We show that in this context it implies the first half of the Arithmetic Mirror Symmetry Conjecture of Gritsenko and Nikulin \cite{gritsnik}.

We put things together in Section $6$, where we find that our modifications of the 
Baily-Borel compactification form a big class that can be described in a uniform manner
and has the Baily-Borel compactification and Mumford's toroidal compactifications
as extremal cases. We remark here that this class of compactifications generalizes 
to all locally symmetric varieties, just like the Baily-Borel and  toroidal 
compactifications, but since we presently see no applications  of this generalization 
we have not endeavoured to set up things in this generality.

An arithmetic arrangement gives rise to a compactification of this type
and it is this case that is discussed in Section $7$. It is here where our main results 
are to be found.  

A new technique must, if not guided, at least be motivated by examples 
and the one developed here is no exception.  This is why
we have included three sections with worked examples to see the theory in action. (Yet 
they are all of the same type, namely moduli spaces of $K3$-surfaces $S$ equipped
with an embedding of a certain lattice in $\pic (S)$, subject to some geometric
condition which is generically satisfied.) In Section $7$
we do this for the geometric invariant theory of $K3$-surfaces of small degree. The
GIT quotient  of the space of plane sextics resp.\ quartic surfaces
(yielding compactifications of the moduli space of nonunigonal $K3$-surfaces of 
degree $2$ and the nonhyperelliptic $K3$-surfaces of degree $4$) has been worked out by 
Shah \cite{shah:deg2}, \cite{shah:deg4} and we recover this quotient as a stratified space in the above setting. (We do not identify here  GIT strata with our strata, but we intend
to do that in a future paper, in which we shall discuss these examples in more detail.) 
As far as we know, no one has done this for the
space of complete intersections of a quadric and a cubic in $\PP^4$ resp.\ of three
quadrics in $\PP^5$ (which give  compactifications of the moduli space of 
nonhyperelliptic $K3$-surfaces of degree $6$ and $8$ respectively). Yet our technique  
predicts what the GIT boundary will be as a stratified space (their dimensions will 
be $3$ and $2$ respectively). 
Section $8$ treats the moduli spaces of `digonal' $K3$-surfaces and of 
Enriques surfaces in this spirit. The last case had already been treated by 
Sterk \cite{sterk} in his Ph.D.\ thesis, so here we limit ourself to describing the 
connection
with his work.  Section $10$ does not involve geometric invariant theory, but
concerns an application to singularity theory:  we obtain a precise 
description of the normalization of the union of the smoothing components of a triangle 
singularity. This, incidentally, is the example that guided us in the first place a
long time ago.

The dependence of this paper on its predecessor \cite{I} is modest. Strictly speaking
this is only so in Section \ref{section:arriv}, where  we need the notion of a linearized arrangement, which is explained at the beginning of that paper. 

\smallskip
\emph{Acknowledgements} This work has developed over a number of
years---I dare not 
say how many, but the reader can guess---and during this period I have benefited 
from conversations and correspondence with many colleagues, among them Egbert Brieskorn, 
Igor Dolgachev, Gert Heckman, David Mumford, Henry Pinkham, Kyoji Saito, 
Peter Slodowy, Hans Sterk and Jonathan Wahl. I thank them all. 
I thank Igor Dolgachev, Gert Heckman and Hans Sterk once more for comments on a first
version. 

The finishing stage of this work was carried out while I was a visitor of 
the Mathematical Sciences Research Institute at Berkeley and in this capacity 
supported by an NSF grant. I extend my thanks to both the MSRI and the NSF.  

\section*{List of Notation}

\noindent 
Notation of a very local nature is omitted here.

\begin{itemize}
\item[$A$] In Sections \ref{section:tube} and \ref{section:arr} an affine 
space over $L$ efines over $\RR$.
\item[$\bb$] Denotes the Baily-Borel compactification.
\item[$C$] In Section \ref{section:tube} a quadratic cone.
\item[$\Ccal$] A tube domain open in $A$ (in Section \ref{section:tube}).
\item[$C_+$] The convex hull of the set of rational points in the closure of
$C$.
\item[$C_I$] The cone naturally associated to the real isotropic subspace 
$I$ of $V$: a half line in $\wedge^I$ resp.\ a quadratic cone in 
$I\otimes I^\perp/I$ resp.\ $\{ 0\}$ if $I$ is of dimension $2$, $1$ 
or $0$.
\item[$C(\DD)$] The conical locus defined by $\DD$ in $\so (\phi )$, see Definition \ref{def:conical}.
\item[$\G$] An arithmetic group, usually of $G(\QQ)$, the exception being 
Sections \ref{section:tube} and \ref{section:arr}, where $\G$ is an arithmetic subgroup of hyperbolic type.
\item[$\tilde\G$] An arithmetic group of automorphisms of the tube domain $\Ccal$.
\item[$\G_I, \G^I, \G (I)$] The $\G$-counterparts of $G_I,G^I, G (I)$ that are arithmetic in these groups; here $L$ is a $\QQ$-subspace of $V$. 
\item[$\G_\sigma,\G^\sigma$] $\G_\sigma$ is the group of $\g\in\G$ which preserve $\sigma$ 
and $\G^\sigma$ the group of $\g\in\G_\sigma$ which act trivially on $V/V_\sigma$ (see 
\ref{not:sigma}). Here $\G$ is an arithmetic group acting on $\DD$ and $\sigma$ is a 
member of an admissible decomposition of the conical locus of $\DD$.   
\item[$\DD$] A domain of type IV in $\PP (V)$ defined by the form $\phi$. 
\item[$\eta$] Usually an ample line bundle.
\item[$\phi$] A nondegenerate symmetric bilinear form on the vector space $V$ 
defined over $\RR$ or $\QQ$ and of signature $(2,n)$.
\item[$\phi_J$] The form $\phi$ induces in $J^\perp/J$, where $J$ is an 
istropic subspace of $V$.
\item[$G$] The stabilizer of $\DD$ in $\Orth (\phi)(\RR)$.
\item[GIT] abbreviation of Geometric Invariant Theory.
\item[$G_I, G^I$] The group of $g\in G$ that preserve resp.\ act as the 
identity on an istropic subspace $I$ of $V$.
\item[$G(I)$] The subgroup of $\GL (I)$ defined by $G_I/G^I$. 
\item[$G^I$] The stabilizer of the $\RR$-isotropic subspace $I$ of $V$ in $G$.
\item[$H$] Usually a hyperplane. When used as a subscript (like in $\LL_H$
and $\BB_H$) it means that we take the intersection with $H$ (or $\PP(H)$).
It is also used for orbit spaces of such intersections relative to an action
of an arithmetic group (like in $X_H$).
\item[$\Hcal$] An arithmetic arrangement in $V$. 
\item[$I, J$] Usually an isotropic subspace in the inner product space $V$.
\item[$I(\sigma)$] The isotropic center of $\sigma$, see Definition \ref{def:semitoricsystem}. 
\item[$L$] In Sections \ref{section:tube} and \ref{section:arr} a vector space defined over $\QQ$,
equipped with a nondegenerate quadratic form of hyperbolic signature.
\item[$L_\sigma$] The $\Sigma$-support space of a member $\sigma$ of the
decomposition $\Sigma$, see \ref{lemma:lsigma}.
\item[$\LL$] The automorphic line bundle over the type IV domain $\DD$.
\item[$\LL^\times$] The automorphic $\CC^\times$-bundle over the type IV domain $\DD$
(contained in $V$).
\item[$\Lcal$] The automorphic (orbi)line bundle over the Baily-Borel
compactification $X$. Also used for its pull-back on other extensions
of $X$ that dominate $X^\bb$.  
\item[$\Lcal (\Hcal)$] and $\Lcal (\Hcal)^{(l)}$ are coherent sheaves of rank
one defined just before Theorem \ref{thm:contract}.
\item[$\ell$] A line bundle over the abelian torsor $Z(J)$ whose total space contains $Z^\tor$.
\item[$\Lambda$] The $K3$-lattice.
\item[$\Lambda_g$] The orthogonal complement of a primitive vector of 
$\Lambda$ with selfproduct $2g-2$.
\item[$N_I$] The unipotent radical of $G_I$. This is also the kernel of
the homomorphism $G_I\to \Orth (I^\perp/I)\times \GL (I)$.
\item[$N_\sigma$] See the notational convention \ref{not:sigma}.
\item[$\Orth (\phi)$] The orthogonal group defined by the symmetric bilinear form $\phi$.
\item[$P_k$] The linear system of degree $k$ divisors on the projective line $P$.
\item[$P'_k$] The subset $P_k$ defined by the reduced divisors.
\item[$\PP (\; )$] The projective space of one-dimensional subspaces of the 
argument (a vector space).
\item[$\PO (\Hcal )$] The partially ordered set of irreducible
components of intersections of members of the arrangement $\Hcal$.
\item[$\pi_L$] Projection along a subspace $L$ of the vector space $V$: 
$\pi_L: V\to V/L$. It is also used for its projectivization $\PP(V)-\PP(L)\to \PP(V/L)$. 
\item[${}^\ss$] Superscript standing for \emph{the set of stable points of}.
\item[${}^\st$] Superscript standing for \emph{the set of semistable points of}.
\item[$Star$] Given a stratification of a topological space of which $S$ is a member, then $Star(S)$ is the union of the members of that partition having $S$ in their closure. 
\item[$\Sigma$] A decomposition into locally rational cones of
$C_+$ (in Sections \ref{section:tube} and  \ref{section:arr}) or
of the conical locus $C(\DD)$ of $\DD$. When used as a superscript: the associated 
semitoric extension. 
\item[$\Sigma (\Hcal)$] The decomposition defined by the arrangement $\Hcal$.
\item[$V$] A complex vector space with $\QQ$-structure endowed with
a symmetric bilinear form $\phi$ also defined over $\QQ$ and of signature $(2,n)$.
\item[$V_\sigma$] The $\Sigma$-support space of a member $\sigma$ of
an admissible decomposition of the conical locus, see Definition \ref{def:semitoricsystem}. 
\item[$X$] The orbifold $\G\bs\DD$.
\item[$\widetilde{X^\circ}$] Denotes an arrangement blowup, see the discussion following (\ref{prop:merforms}).
\item[$\widehat{X^\circ}$] Denotes the contraction
of an arrangement blowup defined by Theorem \ref{thm:contract}.
\item[$X(L)$] The variety $\G_L\bs\pi_L\DD$.
\item[$Z$] In Section \ref{section:onebc} synonymous with $\G^J\bs \DD$, where
$J$ is a $\QQ$-isotropic plane, in Sections \ref{section:tube} and \ref{section:arr}
essentially of the same form, where $I$ is a $\QQ$-isotropic line.
\item[$Z^\tor$] A toric extension of $Z$.
\item[$Z^L$] In Section \ref{section:onebc} an extension of $Z=\G^J\bs\DD$ defined by a 
$\QQ$-subspace of $J^\perp$ which contains $J$. The boundary is $Z (L)$.
\item[$Z(L)$] In Sections \ref{section:onebc}, \ref{section:tube},
\ref{section:arr} synonymous with $\G^J\bs \pi_L\DD$,
elsewhere, in case the argument is a group, it has the traditional meaning
of the center of that group. 
\item[$\bss$] Formation of a categorical quotient.
\item[${}^\circ$] Superscript (as in $X^\circ$) indicates that we
take the complement of an arrangement in the space in question.
\end{itemize}

\tableofcontents

\section{Structure associated to boundary components}\label{section:bc}
Throughout this paper, $n$ is a nonnegative integer, $V$ a complex
vector space of dimension $n+2$ and $\phi : V\times V\to \CC$ a
nondegenerate symmetric bilinear form. We assume that $(V,\phi)$ has
been defined over $\QQ$ and that $\phi$ has signature $(2,n)$ relative
to this $\QQ$-structure. The corresponding orthogonal group $\Orth
(\phi)$ is therefore an algebraic group defined over $\QQ$. Its Lie
algebra $\so (\phi)$ is the space of endomorphisms of $V$ that are
antisymmetric relative to $\phi$ and we shall often identify the
underlying vector space via $\phi$ with $\wedge^2V$ in the obvious manner.

The subset of $\PP (V)$ defined by $\phi (z,z)=0$ and $\phi
(z,\bar{z})>0$ has two connected components that are interchanged by
complex conjugation. Denote by $G$ the subgroup of $\Orth (\phi)(\RR)$
(of index two) which respects the components. Its Lie algebra is of 
course $\so (\phi)(\RR)$. We choose one of the connected components and 
denote it by $\DD$\footnote{For many applications it is better 
to avoid this choice and work with the union of these components
as a $\Orth (V)(\RR)$-manifold throughout. But for the
development of compactification techniques it is more convenient to
work with a single component.}. It is a Hermitian symmetric domain
for $G$.

Let $\LL^\times$ be the corresponding subset of $V-\{ 0\}$ so
that the projection $\LL^\times\to \DD$ is a $\CC^\times$-bundle. For
every integer $k\in\ZZ$, the character $z\in\CC^\times \mapsto z^k\in
\CC^\times$ defines an equivariant line bundle $\LL^k\to\DD$. We get
$\LL:=\LL^1$  by filling in the zero section of $\LL^\times$ and we
shall refer to this bundle as the \emph{natural automorphic line
bundle} over $\DD$. Filling in the section at infinity yields its dual
$\LL^{-1}$. So a function $f:\LL^\times \to \CC$ that is homogeneous of
degree $-k$ may be considered as section of $\LL^k$. The canonical
bundle of $\DD$ can be identified with $\LL^n$: if $p\in\DD$ is
represented by the line $L\subset V$, then the tangent space $T_p\DD$
is naturally isomorphic to $\Hom (L,L^\perp /L)$. Since $\wedge^n
(L^\perp /L)\cong\wedge^{n+2}V$ canonically, we find
that $\wedge^nT_p^*\DD\cong L^n\otimes\wedge^{n+2}V^*$.

We denote by $\G$ an arithmetic subgroup of $G$, that is, a subgroup of
$G$ that is of finite index in the $G$-stabilizer of some lattice in
$V(\QQ)$. It contains a subgroup of finite index $\G'$  that is neat. 
(This means that the multiplicative subgroup of $\CC^\times$ generated by the
eigenvalues of its elements  has no torsion.) Since $\G'$ acts properly
discretely on $\DD$ and without fixed points, the orbit space $\G'\bs
\DD$ exists as an analytic manifold. We regard the $\G$-quotient of
$\DD$ as an orbifold and denote it by
$X$. An important feature of this structure is that it remembers 
the isotropy groups of orbits. In particular,
it detects the center of $\G$ as the kernel of the action of $\G$ on $\DD$.
The $\G$-quotient of $\LL$ defines a line bundle $\Lcal$ on the orbifold $X$
and $\Lcal^{\otimes n}$ is just the $\G$-quotient of $\LL^n$.

\smallskip
A real isotropic subspace of $V$ is of dimension $\le 2$.
Let us first consider the case of dimension $2$.

\subsection{Isotropic planes} Let $J\subset V$ be a real isotropic
plane. Then $J^\perp/J$ is negative definite of dimension $n$. Any
element $g$ of the $G$-stabilizer $G_J$ of $J$ induces a linear
transformation in $J$ and an orthogonal transformation in $J^\perp/J$.
We thus get a group homomorphism $G_J\to \GL (J)(\RR)\times
\Orth(J^\perp/J)(\RR)$. This homomorphism has  an open subgroup as its
image, whereas its kernel $N_J$ consists of transformations that act
trivially on the successive quotients of the flag $0\subset J\subset
J^\perp\subset V$. So $N_J$ is the unipotent radical of $G_J$ and the
Levi quotient $G_J/N_J$ can be identified with an open subgroup of $\GL
(J)(\RR)\times \Orth(J^\perp/J)(\RR)$.

For $g\in N_J$, the restriction of $g-1$ to $J^\perp$ drops to a linear
map $J^\perp/J\to J$. We thus get a map \[ N_J\to \Hom (J^\perp/J,
J)(\RR)\cong (J\otimes J^\perp/J)(\RR), \] where the last isomorphism
used the nondegenerate form on $J^\perp/J$. This is a surjective
homomorphism to a vector group. Its kernel $Z(N_J)$ consists of the
elements that are the identity on $J^\perp$ and every such element is
of the form 
\[ 
z\in V\mapsto z+ \phi (z,e)f-\phi (z,f)e 
\] for certain $e,f\in J(\RR)$. As this transformation only depends
on the image of $e\wedge f$ in $\wedge^2J$, we thus obtain an identification 
\[ 
\exp :\wedge^2J(\RR)\cong  Z(N_J). 
\] 
Here we identify $\lie (G)=\so (\phi)(\RR)$ with 
$\wedge^2V(\RR)$ so that $\wedge^2J(\RR)$ may be regarded
as an abelian subalgebra of this Lie algebra. The commutator must define a
linear map \[ \wedge^2 (J\otimes J^\perp/J)\to \wedge^2 J. \] This map
turns out to be the obvious one, namely contraction by means of $\phi$.
It is in particular  nondegenerate so that  $Z(N_J)$ is the
center of $N_J$ as the notation suggested. This shows that $N_J$ is a
real Heisenberg group. We shall write $\bar N_J$ for the abelian quotient
$N_J/Z(N_J)$.

The action of $N_J$ on $\DD$ realizes the latter as a Siegel domain of
the third kind. First notice that the choice of the component $\DD$
determines an orientation of $J(\RR)$. This means that the punctured
line $\wedge^2 J(\RR)-\{ 0\}$ has a distinguished component. We shall
denote this component by $C_J$. The corresponding half-plane in
$\wedge^2 J$, $\Ccal_J:=\wedge^2 J(\RR)+\sqrt{-1}C_J\subset \wedge^2
J$, is thought of as a subset of the complex Lie algebra $\so (\phi)$.
Notice that $\exp (\Ccal_J)$ is then a semigroup in the
complexification of $Z(N_J)$ (acting by the same formula on $V$) which
preserves $\DD$.

The set of oriented bases of $J^*(\RR)$ defines a half-sphere in the
Riemann sphere $\PP (J^*)$. This half-sphere is just $\pi_{J^\perp}\DD$
once we identify $V/J^\perp$ with $J^*$. In the diagram \[ \DD\to
\pi_J\DD\to \pi_{J^\perp}\DD, \] the first projection is a half-plane
bundle (of complex dimension one), each fiber being a free orbit of the
half-plane semigroup $\exp (\Ccal_J)$. So $Z(N_J)=\exp (\wedge^2
J(\RR))$ acts by translation parallel to the boundary. The action of
$N_J$ on $\pi_J\DD$ is via $\bar N_J\cong J\otimes J^\perp/J$ and
$\pi_J\DD \to\pi_{J^\perp}\DD$ is a principal $\bar N_J$-bundle. Each
fiber of $\DD\to\pi_{J^\perp}\DD$ is a free orbit of the semigroup
$\exp(\sqrt{-1}C_J)\cdot N_J$

In case $J$ is defined over $\QQ$, then so is $G_J$ and $\G_J:=\G\cap
G_J$ is arithmetic in $G_J$. In particular, $\G_{N_J}:=\G\cap N_J$ is
cocompact in $N_J$ and the image of $\G_J$ in the unitary group of
$J^\perp/J$ is finite. Let us  now assume that $\G$ is neat.
Then this image is reduced to the
identity element so that $\G_{N_J}$ is equal to the subgroup
$\G^J\subset\G_J$ of $\g\in \G_J$ that act trivially on $J$. The
quotient $\G (J):= \G_J/\G^J$ is arithmetic, when regarded as a
subgroup of $GL(J)$. In particular, $X(J^\perp):=\G (J)\bs\pi_{J^\perp}\DD$
is a modular curve. The property of $\G^J$ being cocompact in $N_J$
means that $\G^J\cap Z(N_J)$ is infinite cyclic and that the image of
$\G^J$ in $\bar N_J$ is a lattice. If we form the diagram of orbit
spaces \[ \G^J\bs \DD\to \G^J\bs \pi_J\DD\to \pi_{J^\perp}\DD, \] then
the first morphism is a punctured disk bundle and the second a
principal (real) torus bundle. The torus in question has $(J\otimes
J^\perp/J)(\RR)$ as universal covering. It inherits a complex structure
from each fiber in such a manner that in the complex-analytic category,
the bundle is isogenous to a repeated fiber product of the tautological
bundle of elliptic curves over the half-plane $\pi_{J^\perp}\DD$. A
similar description applies to the diagram \[ \G_J\bs \DD\to \G_J\bs
\pi_J\DD\to  X(J^\perp), \] the only difference being that everything is now
fibered over the modular curve $X(J^\perp)$.

\subsection{Isotropic lines} Suppose that $I$ is a real isotropic line
in $V$. Then the form on $I^\perp/I$  induced by $\phi$ has signature
$(1,n-1)$. For our purpose it is however better to tensor this vector space 
with the line $I$. Then $\phi$ induces a symmetric bilinear map
\[
(I\otimes I^\perp/I)(\RR)\times (I\otimes I^\perp/I)(\RR)
\to I\otimes I
\]
which is still defined over $\RR$. Since the line  
$(I\otimes I) (\RR)$ has a natural orientation, we can still speak of its 
signature. As this signature remains $(1,n-1$), the set of
positive vectors is the union of a cone and its antipode. But the choice 
of $\DD$ now singles out one of these two cones. We shall denote that 
distinguished cone by $C_I$.

The stabilizer $G_I$ acts on $I^\perp/I\otimes I$ with image $\bar
G_I$ the group of real orthogonal transformations of $I^\perp/I\otimes I$ that
preserve $C_I$. The kernel $N_I$ of this action consists of
the transformations in $G_I$ that act trivially on $I^\perp/I\otimes I$. If
$e\in I(\RR)$ is a real generator, and we use it to identify 
$I^\perp/I\otimes I$ with $I^\perp/I$, then any such element is of the form
\[ 
\psi_{e,f}: z\mapsto z+\phi (z,e)f-\phi (z,f)e-\tfrac{1}{2}\phi
(f,f)\phi (z,e)e 
\] 
for some $f\in I^\perp(\RR)$. Conversely, every
element of this form lies in $N_I$. Since $f$ is unique modulo $I$, we
obtain an identification of $N_I$ with a vector group: 
\[ 
\exp: (I\otimes I^\perp/I)(\RR)\cong N_I, \quad (e\otimes
f)\mapsto\psi_{e,f} 
\] 
so that $N_I$ is commutative. We regard here
$I\otimes I^\perp/I\subset \wedge^2V$ as an abelian Lie subalgebra of
$\so (\phi )$ defined over $\RR$. 

The line $I$ determines a
realization of $\DD$ as a tube domain: the projection $\DD\to \pi_I\DD$
is a $G_I$-equivariant isomorphism and $\pi_I\DD$ can be characterized
as a subset of $\PP (V/I)$ as follows. The space $\PP (V/I)-\PP
(I^\perp/I)$ is an affine space for the vector group $\Hom (V/I^\perp,
I^\perp/I)\cong I\otimes I^\perp/I$ and is as such defined over $\RR$.
If we divide out the affine space $\PP (V/I)-\PP (I^\perp/I)$ by its
real translations, then we obtain a vector space (the origin is the
image of the real part of the affine space); this vector space is of
course identifiable with $(I\otimes I^\perp/I)(\RR)$.  The tube
domain $(I\otimes I^\perp/I)(\RR)+\sqrt{-1}C_I\subset I\otimes
I^\perp/I$ exponentiates to the semigroup $N_I\exp (\sqrt{-1}C_I)$ in
the complexification of $N_I$. This semigroup preserves $\DD$. 
Observe that $\PP (V/I^\perp)$ and hence $\pi_{I^\perp}\DD$ is a
singleton.

If $I$ is defined over $\QQ$, then so is $G_I$ and $\G_I$ is arithmetic
in $G_I$. This implies that $\G_{N_I}$ defines a lattice in
$I^\perp/I\otimes I$ (preserved by $\bar\G_I$, of course) and that
$\bar\G_I$ arithmetic in $\bar G_I$.

A real isotropic plane $J$ in $V$ containing $I$ corresponds to a real
isotropic line $\bar J$ in $I^\perp/I$ and vice versa. In that case,
$\wedge^2J$ can be regarded as a line in $I\otimes I^\perp/I$ and via
this inclusion, $C_J$ is contained in the boundary of $C_I$.

\section{The Baily-Borel compactification}\label{section:bb}

The isotropic subspaces defined over $\QQ$ play a central role in the
construction of the Baily-Borel compactification  of $X$. We
observed that every $\QQ$-isotropic line $I\subset V$ defines a cone
$C_I\subset (I\otimes I^\perp/I)(\RR)\subset\wedge^2V(\RR)\cong\so
(\phi )(\RR)$. Similarly, every $\QQ$-isotropic plane $I\subset V$
defines an open half-line $C_I$ in $\wedge^2I(\RR)\subset \so (\phi)(\RR)$.  
It is clear that
these cones are mutually disjoint and that $C_J$ meets the closure of
$C_I$ if and only if $I\subseteq J$. We also consider $\{ 0\}$ as an
isotropic subspace of $V$ and put $C_{\{ 0\}}:=\{ 0\}$ and $N_{\{
0\}}=\{ 1\}$. Let $\Ical$ denote the collection of all $\QQ$-isotropic
subspaces of $V$, including $\{ 0\}$. A partial order on $\Ical$ is
defined by the incidence relations, \emph{not} of the isotropic
subspaces, but of the corresponding cones. For instance, $\{ 0\}< J <I$
if $J$ is an isotropic plane containing the isotropic line $I$. For
$I\in \Ical$, we put \begin{equation*} V_I:= \begin{cases} I^\perp &
\text{ if } I\not=\{ 0\},\\ \{ 0\} & \text{ if } I=\{ 0\}. \end{cases}
\end{equation*} So if $I,J\in \Ical$, then $J\le I$ if and only if
$V_J\subseteq V_I$.

\begin{definition}\label{def:conical} 
We call the subset $\cup_{I\in\Ical} C_I$ of $\so
(\phi)(\RR)$ (a disjoint union) the \emph{conical locus} of $\DD$. We 
shall denote it by $C(\DD)$. 
\end{definition}

This locus depends on the domain $\DD$ with its `$\QQ$-structure'. If
we replace $\DD$ by it complex conjugate (the other component), then
the conical locus is replaced by its antipode.

We write $\G^I$ for the group of $\g\in\G$ that preserve $I$ and act as
the identity on $V/V_I$. It is clear that $\G^{\{ 0\}}=\{ 1\}$. 
When $I\not=\{ 0\}$, then $V/V_I$ can be
identified with the dual of $I$ and so $\G^I$ is then simply the group
$\g\in\G$ that leave $I$ pointwise fixed and so this agrees with
our earlier notation. The group $\G^I$ is an 
extension of an automorphism group
of $C_I$ by $N_I$ and so is in fact equal to $\G_{N_I}$ when $\dim
I\not=1$. For the same reason we have  $\G^I N_I=N_I$ unless $\dim
I=1$. Notice that $J\le I$ implies that $\G^J\subseteq \G^I$. Now
consider the disjoint unions 
\[ 
\DD^\bb:=\coprod_{I\in \Ical}
\pi_{V_I}\DD,\quad (\LL^\times)^\bb:=\coprod_{I\in \Ical}\pi_{V_I}\LL^\times 
\] 
(so the term indexed by $\{ 0\}\in\Ical$ yields
$\DD$ resp.\ $\LL^\times$). Both unions come with an obvious action of
$G(\QQ)$. The group $\G$ has finitely many orbits in $\Ical$ and so the
$\G$-orbit set of $\DD^\bb$ is the union of $X$ and finitely many
modular curves and singletons. The Baily-Borel theory puts a
$G(\QQ)$-invariant topology on these spaces such that $\G\bs
(\LL^\times)^\bb\to \G\bs\DD^\bb$ becomes a $\CC^\times$-bundle over a
compact Hausdorff space. This topology, which is at first sight perhaps 
somewhat unnatural, is almost forced upon us if we want $\G$-automorphic 
forms to have a continuous extension for it.  Here is the definition.

\begin{definition} 
For a subset $K\subset\DD$ resp.\ $K\subset
\LL^\times$ we define its \emph{$(\G,I)$-saturation in} $\DD^\bb$
resp.\ $(\LL^{\times})^\bb$ by \[ K^\bb (\G, I):= \coprod_{J\in \Ical
,J\leq I} \pi_{V_J}(\G^I N_I\exp (\sqrt{-1}C_I)K). \] Here we note that
$\G^I$ normalizes $C_I$ so that $\G^I N_I \exp (\sqrt{-1}C_I)$ is
a semigroup in $\Orth (\phi)$. If $K$ runs over the open subsets of
$\DD$ resp.\ $\LL^\times$ and $I$ over $\Ical$, then these form the
basis of a topology, called the \emph{Satake topology}.
\end{definition}

These topologies depend of course on the $\QQ$-structure on $V$, but
not on the particular arithmetic group $\G$ in $G(\QQ)$ since it can be
shown that $G(\QQ)$ acts on both spaces as a group of homeomorphisms.
We regard $\DD^\bb$ as a ringed space by equipping it with the sheaf
$\Ocal_{\DD^\bb}$ of complex valued continuous functions on open
subsets that are analytic on each piece $\pi_{V_I}\DD$. We do similarly
for $(\LL^\times)^\bb$. For such ringed spaces we borrow the terminology
that is standard use in the complex-analytic category.
For instance, $\Ocal_{\DD^\bb}$ is called the \emph{structure sheaf}
of $\DD^\bb$, a local section of this sheaf is called an 
\emph{analytic function} and the projection $(\LL^\times)^\bb\to \DD^\bb$ 
is a morphism, more specifically, it is an analytic $\CC^\times$-bundle.

For every $k\in\ZZ$ we denote by $(\LL^k)^\bb$ the line
bundle over $\DD^\bb$ defined by the representation of $\CC^\times$
whose character is the $k$th power. It is still true that
$\LL^\bb=(\LL^1)^\bb$ resp.\ $(\LL^{-1})^\bb$ is 
obtained by filling in the zero section resp.\ the section at infinity.

The first steps of the theory yield that $\G\bs (\LL^\times)^\bb$ is a
locally compact Hausdorff space with proper $\CC^\times$-action whose
orbit space is compact and can be identified with $\G\bs \DD^\bb$. The
next step is to show that we land in the category of normal analytic
spaces if we take our quotients in the category of ringed spaces (so
$\G\bs (\LL^\times)^\bb$ is then equipped with the sheaf of
continuous complex valued functions whose pull-back to $(\LL^\times)^\bb$ is 
piecewise analytic). For us the following geometric definition 
of an automorphic form is convenient.

\begin{definition} 
A \emph{$\G$-automorphic form of degree $k\in\ZZ$}
on $\DD$ is a $\G$-invariant morphism 
$(\LL^\times)^\bb\to\CC$ that is homogeneous of degree $-k$, in other
words, a continuous complex valued function on $(\LL^\times)^\bb$ that
is analytic and homogeneous of degree $-k$ on each piece
$\pi_{V_I}\LL^\times$. 
\end{definition}

Such a form can also be thought of as a section of $(\LL^k)^\bb$. There
is the related notion of a meromorphic automorphic form, whose
definition the reader will be able to guess. A standard procedure
produces plenty of such forms (see Proposition \ref{prop:merforms} 
for a generalization). Baily and Borel show that the 
$\G$-automorphic forms thus produced separate the points of
$(\LL^\times)^\bb$, which establishes the major part of

\begin{theorem}[Baily-Borel \cite{bb}] 
If $\G$ is neat, then the ringed space $\G\bs \LL^\bb$
defines a complex-analytic line bundle $\Lcal$ on the compact
analytic space $X^\bb:=\G\bs \DD^\bb$. This bundle is ample so that
$X^\bb$ is in fact projective and the partition $\{
X(V^I)\}_{I\in\Ical}$  that comes with its definition is a
finite decomposition into subvarieties. The sections of its
$k$th tensor power are precisely the $\G$-automorphic forms of degree $k$.
The requirement that $\G$ be neat may be dropped, provided that
$\Lcal$ is understood as an orbiline bundle.
\end{theorem}

Actually a bit of a shortcut is possible: if one only knows that
$\G\bs\DD^\bb$ is compact, then the fact that the $\G$-automorphic
forms separate the points of $\G\bs\DD^\bb$ can be used to show that it is
an analytic space.

\begin{none}[Baily-Borel extension]\label{bbext}
 Locally, the Baily-Borel
compactification exists on an intermediate level as an extension
in a suitable complex-analytic category. 
We explain. For $I\in\Ical$, we consider the open
neighborhood \[ Star(\pi_{V_I}\DD):=\coprod_{J\le I} \pi_{V_J}\DD \] of
$\pi_{V_I}\DD$ in $\DD^\bb$. There is an evident retraction
$Star(\pi_{V_I}\DD)\to \pi_{V_I}\DD$ (as a morphism of ringed spaces)
which is equivariant relative to the action of the  $\G$-stabilizer
$\G_I$ of $I$. If $\G^I\subset\G_I$ denotes the subgroup that acts as
the identity on $\pi_{V_I}\DD$, then $\G^I\bs Star(\pi_{V_I}\DD)$ is
already a normal analytic space and the projection 
\[ 
\G^I\bs Star(\pi_{V_I}\DD)\to \pi_{V_I}\DD 
\] 
is analytic. The group
$\G^I/\G_I$ acts properly discretely on the base and hence also on the
total space. It follows that we have a morphisms of analytic spaces 
\[
\G\bs \DD^\bb\leftarrow\G_I\bs Star(\pi_{V_I}\DD)\to X(V^I).
 \]
that are the identity on $X(V^I)$. Reduction theory tells us that
the first map is an isomorphism over a neighborhood of $X(V^I)$. In other
words, there exists a neighborhood $U_I$ of $\pi_{V_I}\DD$ in
$Star(\pi_{V_I}\DD)$ such that every $\G$-orbit meets that neighborhood
in a $\G_I$-orbit or not at all, so that $\G_I\bs U_I$ maps
isomorphically onto a neighborhood of $X(I^\perp)$ in $X^\bb$. 
We shall refer to the extensions (of normal analytic spaces) 
$\G^I\bs\DD\subset\G^I\bs
Star(\pi_{V_I}\DD)$ and $\G_I\bs\DD\subset\G_I\bs Star(\pi_{V_I}\DD)$
as \emph{Baily-Borel extensions}. \end{none}

\begin{remark} It can be shown that for every nonzero $\QQ$-isotropic
subspace $I\subset V$, the Baily-Borel extension of $\G^I\bs\DD$ is a
normal Stein space. Since the boundary of this extension is of
dimension $\le 1$, this gives for $n=\dim\DD\ge 3$ the a priori
characterization of the Baily-Borel extension as a \emph{Stein
completion} of $\G^I\bs\DD$. \end{remark}

\section{Modification over a one-dimensional boundary 
component}\label{section:onebc}

In this section we fix a $\QQ$-isotropic plane $J$ and we abbreviate 
\[
Z:= \G^J\bs\DD,\quad Z(J):=\G^J\bs\pi_J\DD. 
\]
Recall that $Z$ is
a punctured disk bundle over $Z(J)$ and that 
$Z(J)$ is a relative abelian
variety over the half-plane $\pi_{J^\perp}\DD$. To be precise, $\G^J$
defines a lattice $M\subset (J\otimes J^\perp/J )(\RR)$ so that if
we ignore the complex structure on $Z(J)$ we get a 
$(M \otimes\RR/\ZZ)$-principal bundle. We also recall that the center of the
Heisenberg group $\G^J$ has a preferred generator so that formation
of the commutator therefore defines a nondegenerate symplectic form
$M\times M\to\ZZ$. This form is a positive multiple 
of the symplectic form that is defined by
\[ 
\begin{CD} 
(J\otimes J^\perp/J)\times (J\otimes J^\perp/J) @>{\phi_J}>> J\times J
\to\wedge^2J\cong\CC , 
\end{CD} 
\]
where $\phi_J$ is quadratic form on $J^\perp/J$ induced by $\phi$.

The Baily-Borel extension $Z^\bb$ adds to $Z$ a copy of
$\pi_{J^\perp}\DD$ and its analytic structure can be understood as
follows. Filling in the zero section of the punctured disk bundle
$Z\to Z(J)$ (so a copy of $Z(J)$) produces a disk bundle.
This disk bundle is a simple example of the class of toric extensions 
considered by Mumford et al.\  in \cite{amrt}, and so we denote that bundle 
$Z^\tor\to Z(J)$. The bundle is contained in the total space of a line bundle $\ell$ over $Z(J)$; we can identify $\ell$ with the normal bundle of the 
zero section. The first Chern 
class of $\ell$ is given by the symplectic form above and from what has 
been remarked there it follows that its Riemann form on $J^\perp/J$ is a
positive multiple of $\phi_J$. Since $\phi_J$ is negative definite, the
classical theory of abelian varieties then guarantees that the zero
section can be analytically contracted in $Z^\tor$ along the projection
$Z(J)\to\pi_{J^\perp}\DD$. The result of this contraction is the
Baily-Borel extension $Z^\bb$ of $Z$.

\begin{none}
In fact, any $\QQ$-subspace $L$ of $J^\perp$ that contains $J$ defines a
partial blowup 
\[ 
Z^L\to Z^\bb 
\] 
of the Baily-Borel extension by contracting
$Z(J)$ along the $L(\RR)$-orbits in $Z(J)$ (instead of
contracting along $Z(J)\to\pi_{J^\perp}\DD$). This can be done
analytically for the same reason as before, but perhaps a clearer picture 
is obtained by dividing
out not by $\G^J$, but by $\G^L$, the subgroup of elements of $\G_J$
that act as the identity on $V/L$ and to define a Baily-Borel extension
of $\G^L\bs\DD$ as follows. First
note that $\G^L\bs \DD$ is a punctured disk bundle over
$\G^L\bs\pi_J\DD$, the latter being a bundle of abelian varieties over
$\pi_L\DD$. Filling in the zero section adds a copy of
$\G^L\bs\pi_J\DD$. We obtain the Baily-Borel extension alluded to by
contracting the zero section along the projection
$\G^L\bs\pi_J\DD\to\pi_L\DD$. The abelian group $\G^J/\G^L$ acts
on this Baily-Borel extension properly discontinuously and if we pass to the
quotient space we get the partial blowup described above.

However for our purpose it is better to work in a Satake setting
from the start. This goes as follows. Let $N_L$ denote the group of 
real orthogonal
transformations of $V$ that leave the flag $\{ 0\}\subset J\subset
L\subset V$ invariant and act as the identity on the successive
quotients. This is a Heisenberg subgroup of $N_J$. We then put a
topology on the disjoint union $\DD^L:=\DD\sqcup \pi_L\DD$:
if we define for a subset $K\subset \DD$ its \emph{$L$-saturation} as
the subset 
\[ 
K^L:= N_L\exp(\sqrt{-1}C_J) K \coprod \pi_L(N_L\exp(\sqrt{-1}C_J) K) 
\] 
of $\DD^L$, then the collection of
open subsets of $\DD$ and their $L$-saturations define a
$N_J$-invariant topology on $\DD^L$. If we give it the structure of a
ringed space in the usual manner, then $\G^L\bs\DD^L\to\pi_L\DD$
is a morphism of normal analytic spaces that exhibits $\G^L\bs\DD^L$ as
a `mapping cylinder' of the relative abelian variety 
$\G^L\bs \pi_J\DD\to \pi_L\DD$. The projection
$\pi_L\DD\to\pi_{J^\perp}\DD$ is, if we ignore the complex structure, 
also a principal bundle of the vector group $(J^\perp/L\otimes J)(\RR)$.
The group $\G^J/\G^L$ may be identified with a lattice in
this vector group and hence 
\[
Z(L):=(\G^J/\G^L)\bs\pi_L\DD\to \pi_{J^\perp}\DD 
\]
is an abelian torsor.
The group $\G^J/\G^L$ also acts on the map $\G^L\bs\DD^L\to
\pi_L\DD$. This action is proper and free (because it is so on the
base), so that we have an analytic retraction 
\[
Z^L\to Z(L).
\]
If $L'$ is a $\QQ$-subspace of $L$ containing $J$, then we have a
natural map $Z^{L'}\to Z^L$ that is the identity on $Z$. It is
clear that this is a proper analytic morphism.
\end{none}

\subsection{Arrangements near the boundary component}
We encounter such a situation if we are given a
$\G^J$-invariant collection $\Hcal$ of $\QQ$-hyperplanes of $V$ that 
contain $J$ in
which $\G^J$ has only finitely many orbits. The subspace $L$ of
$J^\perp$ in question is then $\cap_{H\in\Hcal}H\cap J^\perp$.
For each $H\in\Hcal$, the intersection $\DD_H:=\DD\cap \PP(H)$ is a
domain of the same type as $\DD$ whose image in $Z$
is of the same type as $Z$. This image, which we denote by
$Z_H$, is a smooth hypersurface in $Z$. The inclusion of $Z_H$ in 
$Z$ has a Baily-Borel
extension $Z_H^\bb\to Z^\bb$ (as ringed spaces, so analytic) which is
injective and proper. The image is the closure of $Z_H$ in $Z$, which
is in general only a Weil divisor. The situation is different for its
closure in $Z^L$. First notice that the boundary of $\DD_H$ in $\DD^L$
is the hyperplane section $\pi_L\DD_H$ of $\pi_L\DD$ 
and that the $\DD^L$-closure of $\DD_H$, $\DD_H \cup \pi_L\DD_H$, is 
the preimage of $\pi_L\DD_H$ under the retraction $\DD^L\to\pi_L\DD$.
The image of $\pi_L\DD_H$ in the abelian 
torsor $Z(L)$ is an abelian subtorsor $Z(L)_H$ of
codimension one, and this is also the boundary of $Z_H$ in $Z^L$.
It is clear that the $Z^L$-closure of $Z_H$ is 
the preimage of $Z(L)_H$ under the retraction $Z^L\to Z(L)$.
In particular, this closure is a Cartier divisor. This proves the first
part of the  following lemma.

\begin{lemma}\label{lemma:cartier1} 
The $Z^L$-closure of each hypersurface $Z_H$ is a Cartier
divisor so that the collection of these defines an arrangement
in $Z^L$. For some  positive integer $k$, the fractional ideal 
$\sum_H\Ocal_{Z_L}(kZ_H)$  
is generated by its global sections. These global sections
separate the points of the arrangement complement in $Z^L$ and that
arrangement complement is Stein. 
\end{lemma} 
\begin{proof} 
Assume for a moment that $\Hcal$ is a single  $\G^J$-orbit. So if 
$H\in\Hcal$, then $L=H\cap J^\perp$
and $Z^L\to Z(L)$ is a relative elliptic curve.
Choose $e\in J-\{ 0\}$. 
If $H$ is defined by the linear form $f_0$ on $V$, 
then  $g_0(z):= f_0(z)/\phi (z,e)$ is a
well-defined function on $\DD$, which even factors through $\pi_L\DD$.  
We form 
\[
E^{(k)}:=\sum_{g\in \G^L g_0} g^{-k}
\]
The index set of this series is actually a rank two quotient of $\G^L$:
it is in fact a classical Eisenstein series. For $k>2$
the series represents a meromorphic function on $Z(L)$ 
which has a pole of exact order $k$ along $Z(L)_H$ and nowhere else.  
The algebra generated by these Eisenstein series
(over the algebra of holomorphic functions on the upper half-plane 
$\pi_L\DD$) contains a set of generators of $\Ocal_{Z(L)}(kZ(L)_H)$, 
when $k$ is sufficiently large. By viewing these as sections of
$\Ocal_{Z^\bb}(kZ_H)$, the lemma follows in this case.

The general case is easily deduced from this: for each 
$H\in \Hcal$, we get a morphism $Z^L\to Z^{H\cap J^\perp}$ such that 
the $Z^L$-closure of $Z_H$ is the preimage of the 
$Z^{H\cap J^\perp}$-closure of $Z_H$. This implies the second assertion 
of the lemma. As the morphism $Z^L\to Z^{H\cap J^\perp}$ only depends
on the $\G^J$-orbit of $H$, there is a natural morphism from
$Z^L$ to the fiber product of the blowups $Z^{H\cap J^\perp}\to Z^\bb$,
where $H$ runs over the finite set $\G^J\bs \Hcal$. This morphism 
is finite. The arrangement complement in $Z^L$ is the preimage under 
this finite map  of a fiber product of Stein manifolds over a Stein 
manifold and hence Stein.
\end{proof}

\begin{corollary}\label{cor:blowup1}
Denote by $D\subset Z$ the union of the hypersurfaces $Z_H$.
If $L$ is of dimension $\ge 2$, then $Z^L\to Z^\bb$ is the normalized
blowup of the fractional ideal $\Ocal_{Z^\bb}(D)$ and the arrangement
complement in $Z^L$ is the Stein completion of the arrangement complement 
$Z-D$ (that is, 
$Z^L$ is a Stein space which has the same algebra of holomorphic 
functions as $Z-D$).
\end{corollary}
\begin{proof}
According to Lemma \ref{lemma:cartier1}, the arrangement on $Z^L$ is
a Cartier divisor and its complement in $Z^L$ is Stein. 
A meromorphic function on $Z^\bb$ which is 
regular on $Z-D$ defines a meromorphic function on $Z^L$. The polar 
set of the latter is of pure codimension one
everywhere and lies in the preimage of the arrangement. Hence it lies 
in the arrangement on $Z^L$. So the direct image
of $\Ocal_{Z^L}(kD)$ on $Z^\bb$ is $\Ocal_{Z^\bb}(kD)$.
The corollary follows.
\end{proof} 

The arrangement on $Z^L$ is in fact the preimage of an arrangement
on $Z(L)$ under the retraction $Z^L\to Z(L)$; we denote
its complement in $Z^L$ by $(Z^L)^\circ$. So we can form without
difficulty the arrangement blowup of $Z^L$ as in Part I \cite{I} to obtain
\[
\widetilde{Z^\circ}\to Z^L.
\]
It follows from Lemma \ref{lemma:cartier1} that this is the blowup 
of an ideal:

\begin{proposition}
The morphism $\widetilde{Z^\circ}\to Z^\bb$ is the blowup of
the fractional ideal $\sum_H \Ocal_{Z^\bb}(Z_H)$. 
\end{proposition}

It is worth noting that we can also obtain
$\widetilde{Z^\circ}$ as the $\G^J$-quotient of the arrangement blowup
of $\DD^\Hcal\to \DD^L$ (since this is locally the preimage 
of an ordinary analytic blowup under a retraction, the absence of a 
locally compact setting is of cause of concern). 
The exceptional divisors of this blowup are 
indexed by the collection $\PO (\Hcal)$ of intersections
of members from $\Hcal$. This index set is partially ordered by 
inclusion and the exceptional divisors indexed by a subset of 
$\PO (\Hcal)$ have  nonempty common intersection if and only if this 
subset is linearly ordered.

\subsection{Arrangements that define $\QQ$-Cartier divisors}
The closure of $Z_H$ in $Z^{J}$ is the bundle restriction of
the disk bundle $Z^{J}\to Z(J)$ 
to $Z(J)_H\subset Z(J)$. There is a
$\QQ$-linear form $l_H: J^\perp/J\to\CC$ whose zero set is the
hyperplane $H\cap J^\perp/J$ in $J^\perp/J$ and has the property that
the rank two symplectic form 
\[ 
\begin{CD} (J^\perp/J\otimes J)\times
(J^\perp/J\otimes J) @>{l_H\otimes\one\times l_H\otimes\one}>> J\times
J \to \wedge ^2J\cong \CC 
\end{CD} 
\] 
represents the first Chern class
of $Z^\tor_H\cdot A$. This linear form is clearly unique up to sign and
its square, $l_H^2$, is the Riemann form of the divisor $Z^\tor_H\cdot
A$.

A $\G_J$-invariant function $H\in\Hcal\to m_H\in \QQ$ defines a
Cartier divisor $Z_\Hcal(m)$ on $Z$ (in an orbi-sense) and a Weil divisor on $Z^\bb$
which we also denote by $Z_\Hcal(m)$.

The following proposition appears in a somewhat different guise in a recent paper by Bruinier and Freitag (Proposition 4.4 of \cite{bruinier}). 

\begin{proposition}\label{prop:principal1}
The Weil divisor $Z_\Hcal(m)$ is a $\QQ$-Cartier divisor if and
only if the quadratic form $\sum_{H\in\G^J\bs \Hcal} n_H l_H^2$ on
$J^\perp/J$ (where the sum is over a system of representive
$\G^J$-orbits in $\Hcal$) is proportional to $\phi_J$. 
\end{proposition}
\begin{proof} 
The proof is standard (see also Section 6 of \cite{I}):
if $Z_\Hcal(m)$ is principal on $Z^\bb$, then some multiple must
be defined by a section of a tensor power of $\ell$. Comparing Chern
classes then gives the proportionality assertion. If conversely the
quadratic forms are proportional, then the restriction of the strict
transform of $Z_\Hcal(m)$ in $Z^J$ to $A$ and $\ell$ have
proportional Chern classes. So a nonzero multiple of the former is a
divisor of a tensor power of $\ell$ translated over an element of
$M \otimes \RR/\ZZ$. It is not hard to check that this
translation lies in $M \otimes \QQ/\ZZ$. So a nonzero multiple of
the strict transform of $Z_\Hcal(m)$ on $Z^J$ is the pull-back along
$Z^J\to A$ of a divisor of a tensor power of $\ell$. This means that
this multiple is the divisor of a holomorphic function on $Z^\bb$.
\end{proof}

\begin{corollary}\label{cor:noproduct}
If $\Hcal$ is nonempty and such that the union of the hypersurfaces 
$Z_H$ supports an effective principal divisor, then the common intersection
of $J^\perp$ and the hyperplanes from $\Hcal$ is equal to $J$.
\end{corollary}
\begin{proof} 
This is because a Riemann form of such a divisor will
vanish on $\cap_{H\in\Hcal} H\cap J^\perp/J$.
\end{proof}

\subsection{Application to moduli spaces of $K3$ 
surfaces}\label{example:k3} 
We digress to show that this corollary almost immediately implies that 
the moduli space of polarized  $K3$-surfaces of given genus  and 
\emph{without} rational double points cannot be affine.

It is well-known that for any integer $g\ge 2$ the
moduli space of primitively polarized $K3$ surfaces of genus $g$ 
($=\text{degree}2g-2$) is isomorphic to an
arithmetic quotient of an arrangement complement, the isomorphism being
induced by a period mapping. To be precise, consider the $K3$-lattice
\[ \Lambda := E_8(-1)\perp E_8(-1)\perp U\perp U\perp U. \] Here
$E_8(-1)$ stands for the $E_8$ lattice whose quadratic form has been
multiplied by $-1$ (so this lattice is negative definite of rank $8$)
and $U$ is the hyperbolic lattice of rank two. The latter has isotropic
generators $e,f$ with inner product $1$. Consider the vector
$h_g:=e_3+(g-1)f_3$ in the last hyperbolic summand. It is indivisible
and has self product $2g-2$. Its orthogonal complement is 
\[
\Lambda_g\cong E_8(-1)\perp E_8(-1)\perp U\perp U\perp I(2-2g), 
\]
where $I(2-2g)$ is the rank one lattice with a generator whose self
product is $2-2g$ (in this case $e_3+(1-g)f_3$). Notice that its signature 
is
$(2,19)$. Let $\DD_g\subset\PP (\Lambda_g\otimes\CC)$ be as usual and
let $\G_g$ be the group of $\g\in\Orth (\Lambda)$ that fix $h_g$ and
leave $\DD_g$ invariant. Consider the collection $N_g$ of
$(-2)$-vectors in $\Lambda_g$ which span with $h_g$ a primitive
sublattice of $\Lambda$. (In case $g>2$, all $(-2)$-vectors span with
$h_g$ a primitive sublattice, but for $g=2$ that is not the case: take
for instance $e_3-f_3$.) It is known that $N_g$ is an orbit under the
group $\G_g$. The collection of hyperplanes in $\Lambda_g\otimes \CC$
orthogonal to a member of $N_g$ is an arithmetic arrangement and
determines therefore an irreducible divisor $D_g$ on $X_g:=\G_g\bs
\DD_g$. Its complement $X_g-D_g$ can be identified with the coarse
moduli space of $K3$ surfaces of genus  $g$. The question has been
raised whether $D_g$ is definable by an automorphic form.
Nikulin \cite{nikulin}  proved  by means of a Koecher 
principle that the answer is no for arbitrary large values of $g$. 
We shall see right away the necessary condition stated in Corollary 
\ref{cor:noproduct} fails for all $g$, so
that the answer is no \emph{always}. This necessary condition comes
down to the requirement that for every isotropic rank $2$ sublattice
$\Pi\subset \Lambda^g$, the image of $\Pi^\perp\cap N_g$ in the negative
definite lattice $\Pi^\perp\cap \Lambda_g/\Pi$ spans a sublattice of
maximal rank. This is clearly not the case: for the simplest choice
$\Pi=\ZZ e_1+\ZZ e_2$ we have 
\[ 
\Pi^\perp\cap\Lambda_g/\Pi\cap\Lambda_g\cong 
E_8(-1)\perp E_8(-1)\perp I(2-2g), 
\] 
and the image of $N_g\cap \Pi^\perp$ herein consists of all the
$(-2)$-vectors in the first two summands. So this image spans a
sublattice of corank one.

\section{Semitoric extension of certain tube
domains}\label{section:tube}

We shall describe a class of modifications of the Baily-Borel extension
of tube domain quotients in terms of combinatorial data. This class was
introduced in \cite{looij:newcomp}, \cite{looij:semitoric} and there
baptized the class of semitoric extensions. It includes both
Mumford's toric extensions and the Baily-Borel extension as special cases.

\begin{definition} Let $L$ be a finite dimensional complex vector space
defined over $\QQ$. Recall that a convex cone in  $L(\RR)$ is called
\emph{nondegenerate} if it does not contain an affine line. We shall call
a cone \emph{rational} if it is the convex cone spanned by a \emph{finite} 
subset of $L(\QQ)$. We say that a finite collection $\Sigma$ of nondegenerate
rational convex cones in  $L(\RR)$ is a \emph{rational cone system} if
a face of a member of $\Sigma$ belongs to $\Sigma$ and any two members
of $\Sigma$ meet along a common face. \end{definition}

If $L$ is obtained from tensoring the cocharacter group of an algebraic
torus with $\CC$, then a rational cone system $\Sigma$ in $L(\RR)$
determines a normal torus embedding of that torus. This construction is
well-known, but since we prefer to work with a torsor of such a torus,
or rather with its universal cover, it may be worthwhile to present
this from a corresponding, inherently analytic, point of view. One
reason being that this will make the connection with the Baily-Borel
compactification clearer.

So let us start with an affine space $A$ over $L$ and a lattice
$L (\ZZ)\subset L(\QQ)$. Then $L/L (\ZZ)$ is an algebraic torus and
$L (\ZZ)\bs A$ is a principal homogeneous space for this torus. A
rational cone system $\Sigma$ in $L(\RR)$ determines a normal affine
torus embedding $L (\ZZ)\bs A\subset (L (\ZZ)\bs A)^\Sigma$ as follows.
For $\sigma\in\Sigma$, denote by $\pi_\sigma: A\to \la\sigma\ra_\CC\bs
A$ the projection along the complex span of $\sigma$ and consider the
finite disjoint union 
\[ 
A^\Sigma :=\coprod_{\sigma\in\Sigma} \la\sigma\ra_\CC\bs A . 
\] 
Notice that $L$ acts on $A^\Sigma$. For
every $K\subset A$ and $\sigma\in\Sigma$, let 
\[
K(\sigma):=\coprod_{\tau \le\sigma}\pi_{\tau}
(\la\sigma\ra_\RR+\sqrt{-1}\sigma +K) \subset A^\Sigma. 
\] 
If $K$ runs over the open subsets of $A$ and $\sigma$ over
$\Sigma$, then we get a $L$-invariant topology on $A^\Sigma$. We make
$A^\Sigma$ a ringed space by endowing it with the sheaf of complex
valued continuous functions that are analytic on its affine
pieces.  The lattice $L (\ZZ)$ acts as a group of automorphisms
of this ringed space and if we pass to the orbit space in this
category, then we get a normal analytic variety that can be identified
with the torus embedding $(L (\ZZ)\bs A)^\Sigma$.

We need local versions of the above notions.

\begin{definition} 
Let $C\subset L(\RR )$  be an open nondegenerate
convex cone and denote by $C_+\supset C$ the convex hull of $L(\QQ)\cap
\overline{C}$ (this is a nondegenerate cone also). 
A subset of $C_+$ is said to be a \emph{locally rational
cone in $C_+$} if its intersection with every rational subcone of $C_+$
is rational. A collection of convex cones in $L(\RR )$ is said to be a
\emph{locally rational decomposition of $C_+$} if the union of these
cones is $C_+$ and the restriction to every rational subcone of $C_+$
is a rational cone system. 
\end{definition}

The coarsest locally rational decomposition of $C$ is the collection of
faces of $C_+$ (this includes $C_+$ itself); for reasons that will
become clear below, we shall denote this facial decomposition of $C_+$
by $\Sigma (\bb)$ (the abbreviation of Baily-Borel).

More interesting examples involve a cone that is homogeneous under a
real semisimple group defined over $\QQ$. The following special case is
relevant for what follows.

\begin{example}[See \cite{looij:semitoric}]\label{example:qcone} Assume
$C$ is a quadratic cone. More precisely, assume there exists a
nondegenerate symmetric bilinear form $L\times L\to \CC$,
$(x,y)\mapsto x\cdot y$, defined
over $\QQ$ and of hyperbolic signature $(1,\dim L-1)$ such that $C$ is
a connected component of the set of $x\in L(\RR)$ with $x\cdot x>0$. Let $\G$ be
an arithmetic subgroup of the orthogonal group of $L(\QQ)$ which
preserves $C$. Here are two constructions of $\G$-invariant locally
rational  decomposition of $C_+$.

(a) Given a finite union $\Ocal$ of $\G$-orbits in $L(\QQ)\cap C_+$,
let $\Sigma$ be the coarsest decomposition of $C_+$ which is closed
under `taking faces' and is such that $\inf_{p\in\Ocal}(-\cdot p)$ is
linear on its members. Then $\Sigma$ is a locally rational
decomposition of $C_+$. If $\Ocal$ is a single $\G$-orbit, then the
members of $\Sigma$ are rational polydral cones if and only if
$\Ocal\subset C$. If $\Ocal$ is a regular orbit, then a maximal member
of $\Sigma$ is a fundamental domain for the action of $\G$ on $C_+$.

(b) Let $\Hcal$ be a collection of $\QQ$-hyperplanes of $L$ which meet
$C$ such that the corresponding subset of the Grassmannian of $L$
is a finite union of $\G$-orbits. Then this collection is locally
finite on $C$ and decomposes $C_+$ into locally rational decomposition
$\Sigma (\Hcal)$. \end{example}

In a similar fashion, a locally rational  decomposition of $C_+$ will
define a \emph{semitoric embedding} (of an open subset of this torus).
In its full generality, this notion is a bit involved, but since we only
need it for the case of Example \ref{example:qcone}, we restrict
ourself to that case. So in the remainder of this section (as well in the next)
we are in

\begin{setting}\label{set}
Let $(L, \cdot )$, $C$ and $\G$ be as in  Example \ref{example:qcone}:
we are given a nondegenerate symmetric bilinear form $L\times L\to \CC$,
$(x,y)\mapsto x\cdot y$, defined
over $\QQ$ and of hyperbolic signature $(1,\dim L-1)$, a connected component 
$C$ of the set of $x\in L(\RR)$ with $x\cdot x>0$ and 
an arithmetic subgroup $\G$ of the orthogonal group of $L(\QQ)$ which
preserves $C$. We let $\Sigma$ be a $\G$-invariant locally rational decomposition of $C_+$. 
We further assume given an affine space $A$
over $L$ \emph{defined over $\RR$} and a group $\tilde\G$ of affine-linear
transformations of $A$ whose translation subgroup $L (\ZZ)$ is a
lattice in $L(\QQ)$ with $\G$ as linear quotient. Since $A$ is
defined over $\RR$ we can identify it with $A(\RR)\times
\sqrt{-1}L(\RR)$. Its open subset $\Ccal:=A(\RR)\times\sqrt{-1}C$ is a
tube domain invariant under the action $\tilde\G$. That action is proper
since the action of $\G$ on $C$ is. So the orbit space
$Z:=\tilde\G\bs\Ccal$ is in a natural manner a normal analytic variety. 
\end{setting}

The semitoric embeddings we are about to define are extensions of 
$Z$ as a normal analytic variety.
But before we start, it is useful to understand the
$\G$-stabilizer $\G_\sigma$ of a $\sigma\in\Sigma$, at least up
to a subgroup of finite index. In case $\sigma\in\Sigma$ meets $C$,
then the subgroup of $\G$ that acts as the identity on the orthogonal
complement of $\sigma$ is of finite index in the $\G$-stabilizer
$\G_\sigma$ of $\sigma$ (because that orthogonal complement is
negative definite) and so  $\tilde\G_\sigma$ acts on $\la\sigma\ra_\CC\bs
A$ via the extension of a finite group by the translation group
$L (\ZZ) /L (\ZZ)\cap \la\sigma\ra_\RR$.

In case $\sigma\in\Sigma$ spans an isotropic line $I$ in $L$, then we
have a homomorphism \[ \exp :I\otimes I^\perp/I\to \SO (L) \] which
assigns to $e\otimes f$ the transformation 
\[ 
l\mapsto l+ (l\cdot e)\tilde f-(l\cdot \tilde f)e+\tfrac{1}{2}
(\tilde f\cdot\tilde f)(l\cdot  e)e, 
\] 
where
$\tilde f\in I^\perp$ is a lift of $f$. The image is the unipotent
radical of the $\SO (L)$-stabilizer of $I$. The generator $e_\sigma$ of
the semigroup  $\sigma\cap L (\ZZ)$ enables us to identify $I\otimes
I^\perp/I$ with $I^\perp/I$. The group of $f\in I^\perp/I$ that
correspond via this identification to an element of $\G$ make up a
lattice (contained in the image of $L (\ZZ)\cap I^\perp$) and this
subgroup of $\G$ is of finite index in $\G_\sigma$. Let us, for the
sake of simplicity, assume that we have equality here: so if $\g\in \G$
acts trivially on $I^\perp$ or $L/I$, then $\g =1$. The
$\tilde\G$-stabilizer $\tilde\G_\sigma$ of $\sigma$ is an extension of this
abelian group by $L (\ZZ)$. This extension does not split: it is a
Heisenberg group with center generated by the translation over
$e_\sigma$. Consider now the restriction of the chain of affine maps
$A\to I\bs A\to I^\perp\bs A$ to $\Ccal$:
\[
\Ccal\to\pi_I\Ccal\to\pi_{I^\perp}\Ccal. 
\] 
The last space is a copy of
an upper half-plane on which $\tilde\G_\sigma$  acts via an
infinite cyclic group of translations (the image of $L (\ZZ)$ in
$L/I^\perp$). Dually, the first projection is an upper half-plane
bundle and the group of elements of $\tilde\G_\sigma$ that act as the
identity on $I\bs A$ is generated by the translation $e_\sigma$. The
second projection is a fibration by affine spaces: it is the
restriction of $I\bs A\to I^\perp\bs A$ to $\pi_{I^\perp}\Ccal$. The
Heisenberg group $\tilde\G_\sigma$ acts on each fiber of $I\bs A\to
I^\perp\bs A$ as a lattice of translations (relative to the  underlying
\emph{real} affine space) with the center acting trivially. So if we pass
to orbit spaces with respect to $\tilde\G_\sigma$, we find that
$\tilde\G_\sigma\bs \Ccal$ is a punctured disk bundle over
$\tilde\G_\sigma\bs \pi_I\Ccal$ and that the latter is a bundle of
complex tori (isogenous to a repeated fiber product of the bundle of 
elliptic curves with itself) over
the punctured disk $\tilde\G_\sigma\bs\pi_{I^\perp}\Ccal$.  Compare this
with our discussion of the isotropic planes in Section
\ref{section:bb}.

\begin{lemmadef}\label{lemma:lsigma} 
Let $\sigma\in\Sigma$. If
$\tau\in\Sigma$ is such that $\tau\ge \sigma$ and $\tau\cap
C\not=\emptyset$, then the common zero set of the real linear forms on
$V$ which are nonnegative on $\tau$ and zero on $\sigma$ is a subspace
that is independent of $\tau$. We call this subspace the
\emph{$\Sigma$-support space of $\sigma$} and denote it by $L_\sigma$.
It is the complex span of $\sigma$ unless $\sigma$ spans an isotropic
line $I$, in which case $L_\sigma$ is the intersection of the
hyperplane $I^\perp$ and the complex linear spaces $\la\tau\ra_\CC$
with $\tau\in\Sigma$ and $\tau >\sigma$. 
\end{lemmadef} 
\begin{proof}
These assertions are clear in case $\sigma$ is not an isotropic half
line. So let us assume that $\sigma$ spans an isotropic line $I$. We
observe that the image of $C$ under the projection $\pi_I: L\to L/I$ is
a real half-space bounded by $(I^\perp/I)(\RR)$. The image of $C_+$
under this projection is the union of that half-space plus the origin;
we prefer to think of this image as a cone over a real affine space
$A_I(\RR)$ over $(I\otimes I^\perp/I)(\RR)$. As we have seen, a
subgroup of finite index of the stabilizer $\G_\sigma$ acts on this
affine space as a lattice of translations. Any $\tau\in\Sigma$ with
$\sigma<\tau$ maps under $\pi_I$ to a (possibly degenerate) convex cone
in $\pi_I(C_+)$ and thus determines a rational polyhedron in
$A_I(\RR)$. These cones decompose $A_I(\RR)$ in a locally finite
$\G_\sigma$-invariant manner. A minimal member of this
decomposition must be an affine subspace of $A_I(\RR)$. If $U\subset
I^\perp/I (\RR)$ denotes its translation space, then it clear that this
is contained in the translation space of the members of the
decomposition of $A_I(\RR)$ that are adjacent to this minimal member.
With induction it follows that the whole decomposition is invariant
under $U$ and is the preimage of a decomposition of $U\bs A_I(\RR)$
into bounded rational polyhedra. The assertions now follow easily with
$L_\sigma$ characterized by the property that $L_\sigma /I\otimes I$ is
the complexification of $U$. \end{proof}

So the support spaces of $\Sigma (bb)$ are $V$, $\{ 0\}$ and the
hyperplanes orthogonal to a $\QQ$-isotropic line.

\smallskip We define an extension of $\Ccal$ as a ringed space. Let 
\[
\pi_\sigma:=\pi_{L_\sigma} :A\to L_\sigma\bs A. 
\] be the projection
along $L_\sigma$ and denote by $\tilde\G^\sigma$ be the subgroup of
$\gamma\in\tilde\G$ which preserve $\sigma$ and each fiber of $\pi_\sigma$.
Then $\tilde\G^\sigma\cdot (L_\sigma (\RR) +\sqrt{-1}\sigma)$ is a
semigroup of transformations of $A$ which preserves $\Ccal$. We put a
topology on the disjoint union 
\[ 
\Ccal^\Sigma:=\coprod_{\sigma\in\Sigma} \pi_\sigma\Ccal 
\] 
as follows. If $\Sigma$ is the facial decomposition, then 
we also write $\Ccal^\bb$ instead of $\Ccal^{\Sigma (\bb)}$.

For every $K\subset \Ccal$ and $\sigma\in\Sigma$, let \[ K^\Sigma
(\tilde\G,\sigma ):=\coprod_{\tau \in\Sigma,\tau \le \sigma}
\pi_\tau(\tilde\G^\sigma (L_\sigma (\RR) +\sqrt{-1}\sigma +K)). \] If $K$
is open in $\Ccal$, then so is $\tilde\G^\sigma (L_\sigma (\RR)
+\sqrt{-1}\sigma +K)$. So if we let $K$ run over the open subsets of
$\Ccal$ and $\sigma$ over $\Sigma$, we get a $\tilde\G$-invariant topology
on $\Ccal^\Sigma$ that is invariant under the $\tilde\G$-stabilizer
$\tilde\G_\sigma$ of $\sigma$. We make $\Ccal^\Sigma$ a ringed space by
equipping it with the sheaf of continuous, piecewise analytic complex
valued functions. We denote the $\tilde\G$-orbit space of $\Ccal^\Sigma$
(as a ringed space) by $Z^\Sigma$. It is clear that $Z^\Sigma$
naturally decomposes into finitely many subspaces of the type 
\[
Z(\sigma):=\tilde\G_\sigma\bs\pi_\sigma\Ccal .
\] 
Let us first get a rough  picture  of the incidence relations. For
every $\sigma\in\Sigma$, the set $A(\RR)+\sqrt{-1}(C\cap Star
(\sigma))$ is open in $\Ccal$ and has the virtue that if a $\tilde\G$-orbit
meets $\Ccal_\sigma$, then it meets the latter in a
$\tilde\G_\sigma$-orbit. This is also true for the interior of the
closure of this set in $\Ccal^\Sigma$, 
\[
U_\sigma :=\coprod_{\tau\le\sigma} \pi_\tau
(A(\RR)+\sqrt{-1}( Star (\sigma)\cap C)). 
\] 
This is an open neighborhood of $\pi_\sigma \Ccal$
that comes with a natural retraction on $\pi_\sigma \Ccal$ (a morphism
in the category of ringed spaces). Notice that $U_\sigma$ and
$U_\tau$ are disjoint unless $\sigma \cup \tau$ is contained in
a member of $\Sigma$, in which case their intersection is
$U_{\sigma\vee\tau}$. The map 
\[
\tilde\G_\sigma\bs U_\sigma\to Z^\Sigma 
\] 
is an isomorphism onto an open neighborhood $U_{Z(\sigma)}$ of
$Z(\sigma)$ that retracts naturally onto $Z(\sigma)$. This clearly
helps us in understanding $Z^\Sigma$ as a ringed space.

We now investigate what the strata $Z(\sigma)$ are like. When
$\sigma\in\Sigma$ meets $C$, then by the previous discussion
$Z(\sigma)$ is a finite quotient of the algebraic torus torsor
$(L_\sigma +L (\ZZ))\bs A$. The situation is somewhat more involved in
case $\sigma\in\Sigma$ spans an isotropic line $I$. Then we have a flag
of vector spaces $I\subseteq L_\sigma\subseteq I^\perp$ which yields
the factorization 
\[ 
\Ccal\to\pi_I\Ccal\to\pi_\sigma\Ccal\to\pi_{I^\perp}\Ccal. 
\] 
The preceding discussion
shows that $\tilde\G_\sigma\bs \Ccal$ is a punctured disc
bundle over $\tilde\G_{\sigma}\bs \pi_I\Ccal$ and that the latter is a
bundle of abelian varieties over the punctured disk $\tilde\G_{\sigma}\bs
\pi_{I^\perp}\Ccal$ which contains the stratum $Z(\sigma)$ as an abelian
subbundle.

So every stratum $Z(\sigma)$ of $Z^\Sigma$ is in fact an analytic
orbifold. For $\sigma\not= \{ 0\}$ it is either of toric or of
(relative) abelian type. This is of course precisely the situation 
that we considered in Section \ref{section:onebc}.

\begin{proposition}\label{prop:analyticspace} 
The ringed space
$Z^\Sigma$ is a normal analytic space. Each stratum $Z(\sigma)$
receives the normal analytic structure exhibited above (so
that it becomes locally closed  for the analytic Zariski topology of
$Z^\Sigma$) and the retraction $U_{Z(\sigma)}\to Z(\sigma)$ is
analytic. 
\end{proposition}

\begin{proof} We only give some indications of the proof; for details
we refer to \cite{looij:semitoric} (a key point of the argument also
appears in \cite{looij:genroot}). According to Baily-Borel \cite{bb},
Theorem 92, we must verify the following properties: 
\begin{enumerate}
\item[(i)] The space $Z^\Sigma$ is locally compact Hausdorff.
\item[(ii)] The topological boundary of each stratum $Z(\sigma)$ is a
union of lower dimensional strata. 
\item[(iii)] The space $Z^\Sigma$ is
normal in the topological sense that it has a basis for its topology
whose members meet $Z$ in connected subsets. 
\item[(iv)] For every stratum $Z(\sigma)$, the restriction map from the
$Z(\sigma)$-restriction of $\Ocal_{Z^\Sigma}$ to the structure sheaf of
$Z(\sigma)$ is surjective. 
\item[(v)] Each point of $Z^\Sigma$ has a
neighborhood on which the local sections of $\Ocal_{Z^\Sigma}$ separate
the points. 
\end{enumerate} 
Property (i) can be proved directly as in
\cite{looij:semitoric}: given a compact subset $K$ of $\Ccal$ one
shows, using the fact that $\G$ has a rational cone in $C_+$ as
fundamental domain, that $K^\Sigma (\tilde\G,\{ 0\})$ maps to a compact
subset of in $Z^\Sigma$. It is clear that this implies that every point
of $Z^\Sigma$ has a compact neighborhood.

The fact that $Z(\sigma)$ is a retract of a neighborhood immediately
implies that distinct points of $Z(\sigma)$ have disjoint neighborhoods
in $Z^\Sigma$. The Hausdorff property for other pairs is somewhat more
subtle, but not that hard. We sketch another proof below that it is
perhaps more instructive.

Properties (ii) and (iii) are easy.

Property (iv) is immediate from the fact that the retraction
$U_{Z(\sigma)}\to Z(\sigma)$ is a morphism of ringed spaces.

Property (v) is proved as follows. Fix a $\sigma\in\Sigma$. Let $f:
A\to\CC$ be a real affine-linear form whose linear part is nonnegative
on $\sigma$ and in some point of $C$ (the latter is automatic if
$\sigma$ meets $C$, of course). Denote by $e^f_0$ the restriction of
$\exp (2\pi\sqrt{-1}f)$ to $\Ccal$. Let $\tau\le\sigma$. If the linear
part of $f$ vanishes on $\tau$, then it vanishes on $L_\tau$ by Lemma
\ref{lemma:lsigma} and so $f$ factors over $L_\tau\bs A$. We denote the
corresponding factor of $e^f_\tau$ on $\pi_\tau \Ccal$ by $e^f_\tau$.
Otherwise $e^f_\tau$ stands for the zero function on $\pi_\tau \Ccal$.
The union of the functions just defined is easily seen to be continuous
on the basic open set $U_{\sigma}$. So that function, denoted 
$e^f$,
is a local section of the structure sheaf. Now assume in addition that
the linear part of $f$ is integral on $L (\ZZ)$. It is clear that $e^f$
then factors through $L (\ZZ)\bs U_{\sigma}$. If $\Ocal$ denotes 
the
$\tilde\G_\sigma$-orbit of $e^f$, then a straightforward estimate shows
that $\sum_{f'\in\Ocal} e^{f'}$ converges on $U_\sigma$ 
uniformly
on subsets of the form $K^\Sigma (\tilde\G,\sigma)$ with $K\subset
U_\sigma$ compact. So this defines a section over $U_{Z(\sigma)}$
of the structure sheaf. One verifies that there are sufficiently many
of these functions to separate the points of $U_{Z(\sigma)}$.
\end{proof}

\begin{remark} One can show that each open subset $U_{Z(\Sigma)}$ is a
Stein space. \end{remark}

Suppose we are given another $\tilde\G$-invariant locally rational
decomposition $\Sigma'$ which is refined by $\Sigma$. So every
$\sigma\in\Sigma$ is contained in a member of $\Sigma'$. If
$\sigma'\in\Sigma'$ is the smallest member with this property, then
$L_\sigma\subset L_{\sigma'}$ and hence there is a natural projection
$\pi_\sigma\Ccal\to \pi_{\sigma'}\Ccal$. It follows that we have an
evident map 
\[ 
\Ccal^\Sigma\to \Ccal^{\Sigma'}  
\] 
of which it is straightforward to verify that it is a morphism 
of ringed spaces (i.e., is continuous and takes the local sections  
of the sheaf on the range to local sections of the sheaf on the domain).
In particular, we have always a morphism $\Ccal^\Sigma\to\Ccal^\bb$. So
if we write $Z^\bb$ for $Z^{\Sigma (\bb)}$, then:

\begin{corollary}\label{prop:propermorphism} The projection
$Z^\Sigma\to Z^{\Sigma'}$ is a proper analytic morphism. In particular,
we have a proper analytic morphism $Z^\Sigma\to Z^\bb$. \end{corollary}

\begin{none} We show that an extension of $Z$ as decribed above lies
between the minimal extension $Z^\bb$ and one of toric type. This can
be used to give another proof that $Z^\Sigma$ and $Z^\bb$ are locally
compact Hausdorff (which is perhaps somewhat easier then the one
sketched above because it brings us in an analytic context at an
earlier stage). First assume that $\Sigma$ is a decomposition into
rational cones. This is the case considered by Mumford and his
collaborators \cite{amrt}: $\Sigma$ determines a normal torus embedding
$L (\ZZ)\bs A\subset (L (\ZZ)\bs A)^\Sigma$ and
$L (\ZZ)\bs\Ccal^\Sigma$ can be identified with the interior of the
closure of $L (\ZZ)\bs\Ccal$ in this torus embedding. So then we are in
the nice situation that $L (\ZZ)\bs\Ccal^\Sigma$ has the structure of a
normal analytic space. The group $\G$ acts on this analytic space and
one verifies that the action is properly discrete. Hence its orbit
space is in a natural way a normal analytic space as well. It is not
hard to see that this orbit space can be identified with $Z^\Sigma$.

We have seen that the boundary strata of $Z^{\Sigma}$ come in two
types: the \emph{toric strata} indexed  by $\sigma\in\Sigma$ which meet
$C$ and the \emph{abelian strata} indexed by an isotropic half-line. A
toric stratum is an algebraic torus modulo a finite group and is
relatively compact in $Z^{\Sigma}$; an abelian stratum is a
hypersurface in $Z^{\Sigma}$ that is fibered over a punctured disc. The
morphism to the Baily-Borel extension, $Z^{\Sigma}\to Z^\bb$, contracts
each toric stratum to a point and is on an abelian stratum the collapse
to a punctured disc.

For a general $\Sigma$, we are in between these two extremes. For
instance, we may intersect $\Sigma$ with any $\G$-invariant
decomposition of $C_+$ into rational cones to produce a refinement
$\tilde\Sigma$ as above. We then have a continous map $\pi
:Z_{\tilde\Sigma}\to Z^\Sigma$. It is not difficult to show that
$Z^{\Sigma}$ has the quotient topology with respect to this map. The
fibers are compact and connected. From this it follows that $\pi$ is
proper and that $Z^\Sigma$ is locally compact Hausdorff. The functions
constructed in our first proof show that
$\pi_*\Ocal_{Z_{\tilde\Sigma}}$ separates the points of $Z^\Sigma$. So
this contraction is analytic. \end{none}

\section{Arrangements on tube domains}\label{section:arr}
 We restrict ourself to the
case where the decomposition of $C_+$ arises as in Example
\ref{example:qcone}-b. To be precise, we assume given a finite union
$\tilde\Hcal$  of $\tilde\G$-orbits in the collection of affine 
hyperplanes of 
$A$ which are defined over $\QQ$ and meet $\Ccal$. 
The corresponding collection $\Hcal$ of linear hyperplanes
of $L$ is then as in Example \ref{example:qcone}-b and so defines a
locally rational  decomposition  $\Sigma (\Hcal)$ of $C_+$. Notice
that the $\Sigma (\Hcal)$-support space of a face spanning an
isotropic line $I$ is equal to $I^\perp$ in case no member of
$\Hcal$ contains $I$ and is the common intersection of the $H\in
\Hcal$ that contain $I$ otherwise. 

For $\tilde H\in\tilde\Hcal$, we denote by 
$\Ccal_{\tilde H}$ the hyperplane section $\Ccal\cap {\tilde H}$.
The $\Ccal^{\Sigma (\Hcal)}$-closure of $\Ccal_{\tilde H}$
is  of the same form as $\Ccal^{\Sigma (\Hcal)}$ itself: 
it meets a stratum $\pi_\sigma\Ccal$ if and only if $\sigma\subset H$
and in that case the intersection is $\pi_\sigma (\Ccal_{\tilde H})$. In
other words, we get the extension of $\Ccal_{\tilde H}$ defined by the
arrangement restriction of  $\tilde\Hcal$ to $\tilde H$. This applies 
in fact to arbitrary intersections: if a collection of such 
hypersurfaces has nonempty intersection, then their common 
intersection meets $\Ccal$ and the 
inclusion of the latter in the former is of the same type as 
$\Ccal\subset\Ccal^{\Sigma (\Hcal)}$. 
If we regard the $\Ccal^{\Sigma (\Hcal)}$-closures of 
the hypersurfaces $\Ccal_{\tilde H}$ as an arrangement on 
$\Ccal^{\Sigma (\Hcal)}$, then the retractions $U_\sigma\to
\pi_\sigma\Ccal$ are compatible with this arrangement: the 
restriction of the arrangement to $U_\sigma$ is the preimage
of its restriction to $\pi_\sigma \Ccal$. So 
the arrangement is in normal directions as
if it were one in an analytic manifold. 

The hyperplane section $\Ccal_{\tilde H}$ maps to an irreducible 
hypersurface $Z_{\tilde H}$ in $Z$. This is a 
Cartier divisor which only depends on the image of ${\tilde H}$ in
$\G\bs\tilde\Hcal$ and since the latter set is
finite, so is the collection divisors thus obtained. The union of the
corresponding collection of hypersurfaces will be denoted by $D$.
We have the following analogue of Lemma \ref{lemma:cartier1}.

\begin{lemma}\label{lemma:cartier2} 
The closure of $Z_{\tilde H}$ in $Z^{\Sigma (\Hcal)}$ is a
Cartier divisor. 
\end{lemma} 
\begin{proof} We show that the
closure of $Z_{\tilde H}$ in $Z^{\Sigma (\Hcal)}$ is Cartier in the
neighborhood $U_{Z(\sigma)}$ of $Z(\sigma)$. For this we use the
retraction $U_\sigma\to\pi_\sigma\Ccal$. The members of
$\tilde\Hcal$ that are parallel to $L_\sigma$ define a locally finite
arrangement on $\pi_\sigma\Ccal$ and thus determine an arrangement of
Cartier divisors on $Z(\sigma)$. The preimage of this arrangement
in $U_{Z(\sigma)}$ under the retraction $U_{Z(\sigma)}\to Z(\sigma)$ is
the trace of the collection of closures of the $Z_H$'s on this open
set. The lemma follows. 
\end{proof}

\begin{lemma}\label{lemma:nonself}
After passing  to a subgroup of $\G$ of finite index,
the closure of any hypersurface $Z_{\tilde H}$ in $Z^{\Sigma (\Hcal
)}$  is without selfintersection.  
\end{lemma}
\begin{proof} 
This is similar to that of Lemma 5.1 in Part I \cite{I} and so we omit it.
\end{proof}

\begin{proposition}\label{prop:semitoricsep}
For some $k>0$,  $\sum_{\tilde H} \Ocal_{Z^{\Sigma (\Hcal)}}(k Z_{\tilde
H})$ is generated by its global sections and these global sections 
separate the points of the arrangement complement in $Z^{\Sigma (\Hcal)}$. 
The arrangement complement in $Z^{\Sigma (\Hcal)}$ is Stein.
\end{proposition}
\begin{proof}
Since  $Z^{\Sigma(\Hcal)}$ is normal, we are already satisfied
if we can find a set of global sections $f_0,\dots ,f_N$ of 
$\sum_{\tilde H} \Ocal_{Z^{\Sigma (\Hcal)}}(kZ_{\tilde H})$ which
generate the latter as a sheaf 
such that the associated meromorphic map $Z^{\Sigma(\Hcal)}\to\PP^N\times
Z^\bb$ is regular on the arrangement complement and has
finite fibers there.  By passing to a subgroup of $\tilde\G$ of finite
index we may and will assume that each $Z_{\tilde H}$ is without 
selfintersection. Since $Z^\bb$ is Stein, its global holomorphic
functions separate its points. 

Fix a $\tilde H_0\in\tilde\Hcal$ and choose an 
affine-linear form $f_0 :A\to\CC$ defined over $\RR$ which has $\tilde
H_0$ as zero hyperplane.  
 
If $K\subset \Ccal$ is 
compact, then a standard estimate shows that the
number of $f\in\tilde\G f_0$ with $N-1 <\sup_K |df| \le N$ is bounded
by a polynomial in $N$ of degree $2\dim L -1$. This implies that for
$k> 2\dim L$, the series  
\[ 
S^{(k)}_{\tilde\G f_0}:= \sum_{f\in\tilde\G f_0} f^{-k} 
\] 
normally converges to a $\tilde\G$-invariant analytic function
on $\Ccal-\cup_{\tilde H\in\tilde\G \{H_0\}}\Ccal_H$.
We claim however that such an estimate holds, not just on $K$, but also
on $\tilde\G (K+ L(\RR)+\sqrt{-1}C_+)$. To see this, we choose 
a rational cone $\Pi\subset C_+$ such that $\G\Pi=C_+$ 
(these exist, see \cite{amrt}, Ch.\ II, \S 4.3, Thm.\ 1) and then observe
that what is needed here is to show that 
for some compact $K'\subset C$, the number of 
$\im (f)\in\G \im(f_0)$ with $N-1 < \sup_{\Pi\cap\G.(K'+C_+)} | \im(f)|
\le N$ is bounded by a polynomial of degree $\dim L-1$ in $N$.
According to \cite{amrt},  Ch.\ II, \S 5 , $\Pi\cap\G.(K'+C_+)$ is of the form $K''+\Pi$ 
with $K''\subset\Pi$ compact. So it is enough to verify our
estimates on $K''+\Pi$. But this follows from the fact
that only finitely many $H\in\Hcal$ meet $\Pi\cap C$. 

As a consequence, $S^{(k)}_{\tilde\G f_0}$
represents a $\tilde\G$-invariant analytic function on the arrangement
complement in $\Ccal^{\Sigma (\Hcal)}$ with a pole of exact order $k$ along
the omitted hypersurfaces, the restriction to $\pi_\sigma\Ccal$
being given by the subseries whose terms are constant on the fibers of
$\pi_\sigma$. We show that this yields enough of such functions to 
separate the points of the arrangement complement in $Z^{\Sigma
(\Hcal)}$. In case $\sigma$ is an isotropic half-line this has been 
established in Section \ref{section:onebc}. So let us focus on the remaining 
case, namely when $\sigma$ meets $C$.
Then $\pi_\sigma\Ccal$ is an affine space
and $\G (\sigma)$ contains the image of $L (\ZZ)$ in $L/L_\sigma$ 
as a subgroup of finite index. Since $L_\sigma$ is an intersection
of members of $\Hcal$, the issue is easily reduced
to the one-dimensional case, as phrased in Lemma \ref{lemma:ratfie2} 
below. 
\end{proof}

\begin{lemma}\label{lemma:ratfie2}
Let $a$ be a positive integer. Then for $k>1$ the series 
\[
S^{(k)}(z):=\sum_{n\in\ZZ} (az+n)^{-k}
\]
is a rational function in $\exp (2\pi\sqrt{-1}z)$. 
If we use the latter to identify $\CC /\ZZ$ with
$\CC^\times$ and if the set of $a$th roots of unity is
regarded as a (reduced) divisor $\mu_a$ on $\CC^\times$, then 
$k(\mu_a)$ is the polar divisor of $S^{(k)}$. 
\end{lemma} 
\begin{proof} 
Normal convergence away from the set of poles is easy.
So we may think of $S^{(k)}$ as a meromorphic function on $\CC^\times$.
That the poles are as decribed is also clear. It remains to see
that $S^{(k)}$ is rational as a function on $\CC^\times$. This
will follow if we show that it is also meromorphic at $0$ and $\infty$.
But since the series defining $S^{(k)}$ converges absolutely and 
uniformly on $|\im (z)|>1$,  $|\re (z)|\le 1$, we see that the latter is
even holomorphic at these points. 
\end{proof}

\begin{corollary}\label{cor:blowup2}
Denote by $D$ the union of the hypersurfaces $Z_H$ in $Z$. If no stratum 
in $Z^{\Sigma (\Hcal)}$ is of codimension one
(which in case $\dim (L)>2$ amounts to: any intersection of members of $\Hcal$ which meets $C_+-\{ 0\}$ is of dimension $\ge 2$),
then $Z^{\Sigma (\Hcal)}$ is the blowup of the fractional ideal
$\Ocal_{Z^\bb}(D)$ and the arrangement complement in 
$Z^{\Sigma (\Hcal)}$ is the Stein completion of $Z-D$.
\end{corollary}
\begin{proof}
According to Lemma \ref{lemma:cartier2}, the arrangement on 
$Z^{\Sigma (\Hcal)}$ is
a Cartier divisor. A meromorphic function on $Z^\bb$ which is 
regular on $Z$ defines a meromorphic function on $Z_{\Hcal
(\Sigma)}$. The polar set of the latter is of pure codimension one
everywhere and lies in the preimage of the arrangement. Hence it lies 
in the arrangement on $Z_{\Hcal(\Sigma)}$. So the direct image
of $\Ocal_{Z^{\Sigma (H)}}(kD)$ on $Z^\bb$ is $\Ocal_{Z^\bb}(kD)$.
The assertions then follow from Proposition \ref{prop:semitoricsep}.
\end{proof}

If we are in the situation of the conclusion of Lemma 
\ref{lemma:nonself}, then we may form the arrangement blowup
of $Z^{\Sigma (\Hcal )}$ relative to the closures of the 
hypersurfaces $\{ Z_{\tilde H}\}_{\tilde
H\in\tilde\Hcal}$. But it is better to work independently of 
$\tilde\G$ and to introduce this blowup
as the $\tilde\G$-quotient of an arrangement blowup 
of $\Ccal^{\Sigma (\Hcal)}$.  
We can perform that blowup as in Part I \cite{I} to obtain
\[
\widetilde{\Ccal^\circ}\to\Ccal^{\Sigma (\Hcal)}.
\]
Since this is locally the preimage of an ordinary analytic blowup
under a retraction, the absence of a locally compact setting
is of no concern to us. The exceptional divisors of this blowup are 
indexed by the collection $\PO (\tilde\Hcal |_\Ccal)$ of intersections
of members from $\tilde\Hcal$ that have nonempty intersection with
$\Ccal$. This index set is partially ordered by inclusion and the 
exceptional divisors indexed by a subset of $\PO (\tilde\Hcal |_\Ccal)$
have  nonempty common intersection if and only if this subset is 
linearly ordered. If we pass to $\tilde\G$-orbit spaces we get an analytic morphism 
\[
\widetilde{Z^\circ}:=\tilde\G\bs\widetilde{\Ccal^\circ}\to Z^{\Sigma (\Hcal )}.
\]
It follows from
the preceding that this is the blowup of an ideal on $Z^\bb$:

\begin{proposition}
The morphism $\widetilde{Z^\circ}\to Z^\bb$ is the blowup of
the fractional ideal $\sum_{\tilde H} \Ocal_{Z^\bb}(Z_{\tilde H})$. 
\end{proposition}

\subsection{Arrangements definable by a product expansion} 
This subsection is not indispensable for the rest of this paper.
We here investigate the situation where the
arrangement complement $Z-D$ is Stein. Suppose we
are given a $\tilde\G$-invariant function ${\tilde H}\in\tilde\Hcal\mapsto 
m_{\tilde H}\in\ZZ$.
This defines a Cartier divisor 
\[
 Z_{\tilde\Hcal} (m):=\sum_{{\tilde H}\in\G\bs\tilde\Hcal}  m_{\tilde H} Z_{\tilde H} 
\] 
on $Z$ which we also regard as a  Weil divisor on
$Z^\bb$. We ask a similar question as in Section \ref{section:onebc}: 
when is $Z_{\tilde\Hcal} (m)$ Cartier on $Z^\bb$? 
Proposition \ref{prop:principal1} gives us immediately a rather
strong necessary condition: recall that $Z^\bb-Z$ has strata
$Z(\sigma)$, with $\sigma$ either a $\QQ$-isotropic half-line of $C_+$
(these have dimension one) and or $\sigma=\{ 0\}$ (this is a
singleton). In the former case $\sigma$, the local structure of $Z^\bb$
near $Z(\sigma)$ is not very complicated: as we have seen that it is
basically like the  $\tilde \G_I$-orbit space of $\Ccal \sqcup
\pi_{I^\perp}\Ccal$, where $I$ is the span of $I$. This is the situation
considered in \ref{prop:principal1}. Let us translate the statement of 
that proposition to the present situation. The $\tilde\G$-stabilizer 
$\tilde\G_I$ of $I$ has finitely many orbits in the set $\tilde\Hcal^I$ of 
$I$-invariant members of
$\tilde\Hcal$. Every ${\tilde H}\in \tilde\Hcal^I$ determines a $\QQ$-hyperplane in
$I^\perp/I$ and has associated to it a  $\QQ$-linear form  given up to
sign defining this hyperplane: $\pm l_{\tilde H}:I^\perp/I\to\CC$. Its square
$l_{\tilde H}^2$ depends in a $\tilde\G_I$-invariant manner on ${\tilde H}$ and in order
that $Z_{\tilde\Hcal}(m)$ is principal near $Z(\sigma)$ a necessary 
and sufficient condition is that the quadratic form 
\[
\sum_{{\tilde H}\in \tilde\G_I\bs\tilde\Hcal^I} m_{\tilde H}l_{\tilde H}^2 
\] 
on $I^\perp /I$ is proportional to the
form that the given hyperbolic form induces on this subquotient of $L$. If $\tilde\Hcal^I$ is
nonempty, then this proportionality  factor must be nonzero and so the
intersection of the hyperplanes taken from $\tilde\Hcal^I$ must be reduced 
to $I$. In particular, the $\Sigma (\Hcal)$-support space of $\sigma$
as defined in \ref{lemma:lsigma} is equal to $I$.

\begin{none}
Let us call a member of $\Sigma (\Hcal)$ with nonempty interior
a \emph{chamber} of $\Sigma (\Hcal)$.
For every chamber $\sigma$ of $\Sigma (\Hcal)$ and ${\tilde H}\in \tilde\Hcal$, there
is a unique affine-linear form $f: A\to \CC$ with (i) ${\tilde H}$ as zero set,
(ii) has a linear part that is positive on $\sigma$ and (iii) whose
orbit under the translation lattice $L (\ZZ)$ is $f+\ZZ$. We denote
that function $f_{\tilde H}^\sigma$ and let $e^{f_{\tilde H}^\sigma}$ have the meaning as
before: the restriction of $\exp (2\pi\sqrt{-1}f_{\tilde H}^\sigma)$ to $\Ccal$.
Since $e^{f_{\tilde H}^\sigma}$ is $L (\ZZ)$-invariant, we may also regard it as
a function on $L (\ZZ)\bs\Ccal$. If $\sigma'$ is another chamber then the collection $\tilde\Hcal (\sigma,\sigma')$ of
members of $\tilde\Hcal$ whose translation space separates $\sigma'$ from
$\sigma$ is a union of finitely many $L (\ZZ)$-orbits. This allows us
to define the product 
\[ 
e^m(\sigma,\sigma'):=\prod_{{\tilde H}\in L (\ZZ)\bs\tilde\Hcal
(\sigma,\sigma')} (-e^{f_{\tilde H}^\sigma})^{m_{\tilde H}}. 
\] 
It is easily checked that if
$\sigma''$ is a third member of $\Sigma (\Hcal)$, then 
\[
e^m(\sigma,\sigma'')=e^m(\sigma,\sigma').e^m(\sigma',\sigma''). 
\] 
It is also clear that $e^m$ is $\G$-invariant: for $\g\in\G$ we have
$e^m(\g\sigma,\g\sigma')=\g e^m(\sigma,\sigma')
(=(\g^{-1})^*e(\sigma,\sigma'))$. This implies that 
\[
e^m(\sigma,\g\g'\sigma)=e^m(\sigma,\g\sigma). \g e^m(\sigma,\g'\sigma). 
\] 
In other words, for a fixed $\sigma$, $\g\in\G\mapsto e^m_\sigma (\g):=
e^m(\sigma,\g\sigma)$ is a $1$-cocycle of $\G$ with values in the group
of quasi-characters of $L (\ZZ)\bs A$. Such a cocycle defines an action
of $\G$ on the trivial line bundle $\CC \times L (\ZZ)\bs A\to
L (\ZZ)\bs A$ by letting $\g\in\G$ send $(t,z)$ to $(e^m_\sigma
(\g^{-1})(z)t,\g z)$. The special form of $e^m_\sigma$ makes that the
action even extends to one of the semidirect product of $\G$ with the
torus $V/ L (\ZZ)$. The restriction of this bundle to $L (\ZZ)\bs
\Ccal$ descends to a  orbiline bundle on $Z$ that we will denote by
$\EE^m$. A section of the latter is given by a function $k:
L (\ZZ)\bs \Ccal\to\CC$ with the property that $\g k=e^m_\sigma
(\g).k$ for all $\g\in\G$. The identity $e^m_{\sigma
'}(\g) .e^m_\sigma (\g)^{-1}=\g e^m(\sigma,\sigma').e^m(\sigma,\sigma')^{-1}$
shows that the cohomology class of $e^m_\sigma$ is independent of
$\sigma$ so that the isomorphism class of $\EE^m$ is well-defined. 

A divisor for $\EE^m$ is obtained as follows.
Consider the product expansion 
\[ 
P^m_\sigma:=\prod_{{\tilde H}\in L (\ZZ)\bs \tilde\Hcal} (1-e^{f_{\tilde H}^\sigma})^{m_{\tilde H}}. 
\] 
It represents an analytic function on $  L(\ZZ)\bs\Ccal$ for
if $K\subset C$ is compact, then the series
$\sum_{{\tilde H}\in  L(\ZZ)\bs \tilde\Hcal} |e^{f_{\tilde H}^\sigma}|$ converges uniformly on
$A+\sqrt{-1}K$ and as is well-known, then the same must hold for 
the corresponding product. 
We notice that $P^m_\sigma$ satisfies the functional equation
\[
\g P^m_\sigma =e^m_\sigma (\g)^{-1}.P^m_\sigma  
\]
and so $P^m_\sigma$ represents a section of the dual of $\EE^m$.
The divisor of this section is clearly $-Z_{\tilde\Hcal} (m)$.
\end{none}

\begin{proposition}\label{prop:weylv} 
The line bundle $\EE^m$ is trivial if and
only if there exists a quasicharacter $e^{\rho}$ on $  L(\ZZ)\bs
\Ccal$ such that $\g e^{\rho}=e_\sigma (\g).e^\rho$. This
quasicharacter is  unique up to constant and in that case $e^{-\rho}.
P^m_\sigma$ represents an analytic function on $Z^\bb$ with divisor 
$Z_{\tilde\Hcal}(m)$. 
\end{proposition} 
\begin{proof} 
The line bundle $\EE^m$ is trivial if and only the line bundle
$\CC\times  L(\ZZ)\bs \Ccal\to  L(\ZZ)\bs \Ccal$ with the given
$\G$-action admits a $\G$-invariant section. Such a section 
amounts to an analytic function $k:  L(\ZZ)\bs \Ccal\to\CC^\times$  
satisfying the functional equation $\g k=e^m_\sigma (\g).k$ for all 
$\g\in\G$. If we have
such a $k$, then any nonzero term of its Fourier development yields a
quasicharacter satisfying the same functional equation. If two
quasicharacters have that property, then their quotient is a
$\G$-invariant quasicharacter. But it is clear that such a
quasicharacter has to be constant. 

Given such a $\rho$, then  
$e^{-\rho}.P^m_\sigma$ is $\tilde\G$-invariant and hence represents 
an analytic function on $Z$. If $Z^\bb-Z$ is of codimension $\ge 2$ in the 
normal analytic space $Z^\bb$, then this function extends analytically to 
$Z^\bb$. But this codimension condition only fails in cases where 
there is nothing to prove, namely when $\dim L\le 1$ or $\dim L=2$ and
the form represents zero, for then $\G$ is finite.
\end{proof}

\begin{remark}
The above proposition can also be phrased in terms of toroidal geometry:
for any quasicharacter $e^\rho :  L(\ZZ)\bs\Ccal\to\CC^\times$, 
we define a continuous function 
$r:C_+\to\RR$ which is $\ZZ$-valued on $L(\ZZ)\cap C_+$ and
piecewise affine-linear relative to the 
decomposition $\Sigma(\Hcal)$ by letting it on $\sigma'\in\Sigma (\Hcal)$ 
be given as the linear part of $\rho$ plus the sum of the linear parts
of the functions $m_{\tilde H}f_{\tilde H}$, where ${\tilde H}$ runs over a system of 
representatives of the (finitely many) $L(\ZZ)$-orbits in 
$\tilde\Hcal (\sigma,\sigma')$. 
Since $f_{\tilde H}$ is constant on ${\tilde H}$, the function $r$ is well-defined and 
continuous. (In case $m_{\tilde H}\ge 0$ for all ${\tilde H}$, then this function is also
convex in the sense that the set of $(x,t)\in C_+\times\RR$ with 
$t\ge r(x)$ has that property.)  Then $e^\rho$ satisfies the functional 
equation if and only if $r$ is $\G$-invariant. 
It is not hard to show that a converse also holds: any $\G$-invariant 
$\RR$-valued continuous function $r:C_+\to\RR$ that is
$\ZZ$-valued on $  L(\ZZ)\cap C_+$ and
piecewise affine-linear with respect to the decomposition $\Sigma(\Hcal)$
comes from a function $m: \tilde\G\bs \tilde\Hcal\to\ZZ$ to which the above
proposition applies. 
\end{remark}

Proposition \ref{prop:weylv} yields strong restrictions on the chamber structure: 

\begin{proposition}\label{chamberstructure}
Suppose that we are in the situation where Proposition \ref{prop:weylv} applies:
the closure of $D$ in $Z^\bb$ is the support of
a (not necessarily) effective Cartier divisor.
Assume also that $\dim (L)>2$ and that $\Hcal$ is nonempty.
Given a chamber $\sigma$ of $\Sigma (\Hcal)$, let $v_\sigma\in L(\QQ)$
be characterized by the property that the inner product with $v_\sigma$
is the linear part of $\rho$, where $e^{\rho}$ is as in \ref{prop:weylv}. Then:
\begin{enumerate}
\item[(i)] the $\G$-stabilizers of $\sigma$ and $v_\sigma$ coincide,
\item[(ii)] $\sigma$ is a rational cone if and only if $v_\sigma\cdot v_\sigma >0$, and
\item[(iii)] in case $\sigma$ is not a rational cone, then (a) $v_\sigma$ is isotropic and nonzero, (b)  $v_\sigma$ lies in the span of an isotropic edge of $\sigma$, (c) $\G_\sigma$ contains a finite index subgroup which is free abelian of rank $\dim (L)-2$ and (d) $\sigma$ is a locally polyhedral cone in the half space  that contains $C_+$ and has $v_\sigma$ in its boundary.
\end{enumerate}
\end{proposition}

\begin{proof}
The first assertion follows from the fact that $v_{\g\sigma}=\g (v_\sigma)$
for all $\g\in\G$.

If $\sigma$ is a rational cone, then the $\G$-stabilizer of $v_\sigma$
must be finite. Since $v_\sigma\in L(\QQ)$ and $\dim (L)>2$, this
can only be if $v_\sigma\cdot v_\sigma >0$. 

Assume now that $\sigma$ is not a rational cone.
We first show that there exists a rational cone $\pi\subset\sigma$ 
such that $\G_\sigma .\pi =\sigma$. By reduction 
theory there exists a rational cone $\Pi\subset C_+$ such that
$\G_\sigma .\Pi =C_+$. We know that the  decomposition 
$\Sigma (\Hcal)|_\Pi $ is finite.
For every piece $P$ of this decomposition that lies in $\G.\sigma$ 
we choose a $\g_P\in \G$ such that $\g_P P\subset\sigma$. If we take for
$\pi\subset\sigma$ the cone spanned by the finitely many 
rational cones $\g_P P$ thus found, then $\pi$ is as desired.

Let $\tau$ be any (one-dimensional) isotropic edge of $\pi$. If $\tau$
is an intersection of members of the arrangement, then it is also
an intersection of finitely many supporting hyperplanes of $\sigma$.
If all isotropic edges of $\pi$ are of this form, then the same must be true
for all edges of $\sigma$ and so the projectivization of $\sigma$ is 
a convex locally polyhedral subset of the projective space of $L$. 
But such a subset is in fact a finite polyhedron, and this
contradicts our assumption that $\sigma$ is not a rational cone.
 
So there exists an isotropic edge $\tau$ that is not an intersection of 
members of the arrangement. According to Corollary \ref{cor:noproduct}, 
$\tau$ is then
not contained in \emph{any} member of the arrangement. This implies that $\sigma$ is invariant under $\G_\tau$. The fixed point subspace 
of $\G_\tau$ in $L$ is equal to the span of $\tau$, and so both (iii-a) and
(iii-b) follow. This shows in particular that 
$\tau$ is unique and that $\G_\sigma=\G_\tau$. Assertion (iii-c)
the follows from the fact that $\G_\tau$ contains a free abelian subgroup of finite index of rank $\dim (L)-2$. Now let $\pi'$ be the intersection
of $\sigma$ and the half spaces whose bounding hyperplane supports $\pi$ and contains $\tau$. Clearly, $\pi$ is then a neighborhood of $\tau$ in $\pi'$ 
and $\G_\sigma\pi'=\sigma$. 
Since the projectivization of $\pi'$ is convex and locally polyhedral
in the projective space of $L$, $\pi '$ is a rational cone. This implies
assertion (iii-d).  
\end{proof}

The two preceding propositions are related to work of Gritsenko and Nikulin.
For this we assume that $\tilde\Hcal$ is a 
\emph{reflection arrangement} in the sense that each member is the fixed point hyperplane of a (real) reflection contained in $\tilde\G$. 
Since $\tilde\Hcal$ is $\tilde\G$-invariant, the 
subgroup $\tilde W$ of $\tilde\G$ generated by these reflections is normal in $\tilde\G$. Then the image $W$ of $\tilde W$ in $\G$ is generated by the $\Hcal$-reflections and is normal in $\G$. 
It is well-known that the latter is a Coxeter group which permutes the chambers 
simply transitively.
Now part (i) of their Arithmetic Mirror Symmetry Conjecture 2.2.4 in \cite{gritsnik} follows from assertions (ii) and (iii-d) of Proposition \ref{chamberstructure}. The triple $(V(\QQ),\G,W)$ is almost the same thing as what Gritsenko and Nikulin call a \emph{reflective hyperbolic lattice} (Definition (1.5.1) of \cite{lorenzian}). They observe (Proposition (1.5.2) of \emph{op.\ cit.}) that there are only finitely many isomorphism types of these in dimension $>2$. They can be enumerated, at least in principle, and this leads to a classication of what they call \emph{Lorentzian} Kac-Moody algebras. 
  
We conclude this section by mentioning the following consequence of Proposition \ref{prop:weylv} (the proof of which we omit):

\begin{proposition}\label{prop:weyl}
When the arrangement is reflective in the above sense, then the closure of $D$ in $Z^\bb$ supports an effective Cartier divisor if and only if we can assign in a $\G$-equivariant manner  to each chamber $\sigma$ a nonzero vector 
$v_\sigma\in\sigma\cap L(\QQ)$ which does not lie in any member of $\Hcal$.
\end{proposition}

\begin{example}
A beautiful (and now classical) example to which these results apply is due to Conway (Chapter 27 of \cite{conway}) and Borcherds 
\cite{borcherds:moon} respectively: take for $L(\ZZ)$ any even unimodular lattice of signature $(1,25)$ (there is only one isomorphism class of these), let $\Hcal$ be the collection of all hyperplanes perpendicular to the $(-2)$-vectors in $L(\ZZ)$, $C\subset L(\RR)$ one of the two cones and
let $\G$ be the subgroup of the orthogonal group of $L(\ZZ)$ which preserves
$C$. Since the direct sum of the Leech lattice (with a minus sign) and a hyperbolic lattice is even unimodular of signature $(1,25)$, it follows that there exists a primitive isotropic vector $v\in L(\ZZ)$ in the closure of the 
positive cone such that $v^\perp/\ZZ v$ is isomorphic to (minus) the Leech lattice. It also follows that such vectors make up a $\G$-orbit.
Let us fix such a $v$. Clearly, $v$ is not in any reflection hyperplane. It is easy to see that $\G_v$ maps isomorphically onto the group of affine-linear isometries of $v^\perp/\ZZ v$. One verifies that if $u$ is a
primitive isotropic vector not in $\G v$, then the $(-2)$-vectors in 
$u^\perp/\ZZ u$ span a sublattice of finite index. So if $\sigma$ is the unique chamber that contains $v$, then $v$ is a Weyl vector, by Proposition \ref{prop:weyl}. We can therefore write down a $\G$-invariant infinite product expansion (Borcherds' denominator formula) which defines the discriminant $D$.  
\end{example}

\section{Semi-toric compactification of type IV
domains}\label{section:ivsemitoric}

We return to the situation of Section \ref{section:bb}.

\begin{definition}\label{def:semitoricsystem} We say that a
$\G$-invariant collection $\Sigma$ of cones in the faces of the conical
locus of $\DD$ is an \emph{admissible decomposition of the conical
locus} if \begin{enumerate} \item[(i)] for every $\QQ$-isotropic line
$I$, the members of $\Sigma$ contained in $C_{I,+}$ define a locally
rational decomposition of $C_{I,+}$ and \item[(ii)] if $J$ is a
$\QQ$-isotropic plane (so that the cone $C_{J,+}=\bar C_J$ in
$\wedge^2J$ is member of $\Sigma$) and  $I$ is any $\QQ$-isotropic line
in $J$ (so that $J/I$ is a $\QQ$-isotropic line in $I^\perp/I$), then
the support space of $C_{J,+}$ in $I^\perp/I\otimes I$ relative to the
decomposition $\Sigma|C_{I,+}$ (see \ref{lemma:lsigma}) is independent
of $I$ when we regard that support space as a subspace of $J^\perp$
containing $J$ \end{enumerate} and we then define for any $\sigma\in
\Sigma$ its \emph{$\Sigma$-support space $V_\sigma\subset V$} as
follows: \begin{enumerate} \item[(i)] if $\sigma$ lies in $C_{I,+}$ for
some $\QQ$-isotropic line $I$, but is not an isotropic half-line in
$C_{I,+}$, then $V_\sigma$ is the subspace of $I^\perp$ containing $I$
corresponding to the complex-linear span of $\sigma$, \item[(ii)] if
$\sigma=C_{J,+}$ for some $\QQ$-isotropic plane, then $V_\sigma$ is the
subspace of $J^\perp$ containing $J$ as defined by condition (ii)
above. \end{enumerate} There is a slight ambiguity of notation in case
(ii), for then $V_\sigma$ not only depends on $\sigma$ (or
equivalently, $J$), but also on $\Sigma$. We therefore sometimes write
$J_\Sigma$ for this space instead. The \emph{isotropic center}
$I(\sigma)$ of $\sigma$ is the subspace $V_\sigma\cap V_\sigma^\perp$
of $V_\sigma$. \end{definition}

\begin{example} The simplest example is the
decomposition of the conical locus into its faces. The support spaces
are $\{0\}$ and the orthogonal complements of the nontrivial
$\QQ$-isotropic subspaces. We will see that this decomposition gives rise 
to the Baily-Borel compactification and so we call it the \emph{Baily-Borel
decomposition}.  
\end{example}

\begin{example} Suppose we are given for every $\QQ$-isotropic plane
$J$ a positive rational generator $r_J\in \wedge^2 J (\QQ)$ such that
this assignment is $\G$-equivariant. Fix for the moment a
$\QQ$-isotropic line $I$. For any $\QQ$-isotropic plane $J\supset I$,
$r_J$ can be regarded as an isotropic vector in $(C_I)_+$ via the
inclusion $\wedge^2 J \subset I\otimes I^\perp/I$. The bounded faces of
the convex hull of these isotropic vectors are convex hulls of finite
sets of such vectors. The cones over these define a decomposition
$\Sigma_I$ of $(C_I)_+$ into rational  cones. The union $\Sigma:=\cup_I
\Sigma_I$ is an admissible decomposition of the conical locus into
rational  cones. It has the property that for every $\QQ$-isotropic
plane $J$ we have $J_\Sigma =J$.

Notice that there is a natural choice for the assignment $J\mapsto r_J$
if we are given a $\G$-invariant lattice $V(\ZZ)\subset V(\QQ)$: just
take $r_J$ to be the positive generator of $\wedge^2 J(\ZZ)$.

We can also do a dual construction by taking instead for $\Sigma_I$ the
coarsest decomposition of $C_+$ which is closed under taking faces and
is such that $\inf_{J\supset I}\phi (-,r_J)$ is linear on each member.
Then $\Sigma_I$ is only a locally rational decomposition of $(C_I)_+$.
The union $\Sigma:=\cup_I \Sigma_I$ is admissible as before, and now we
have $J_\Sigma =J^\perp$ for every $\QQ$-isotropic plane $J$.
\end{example}

\begin{example}[Admissible decompositions from arrangements]\label{example:arrdec}
Our main example of interest is the one that arises from an
arithmetic arrangement on $\DD$. First note that every hyperplane
$H\subset V$ defined over $\RR$ of signature $(2,n-1)$ gives a nonempty
hyperplane section $\DD_H :=\PP(H)\cap \DD$ of $\DD$. This is a
symmetric domain for its $G$-stabilizer and a totally geodesic
hypersurface of $\DD$. Now let $H$ run over a finite union $\Hcal$ of
$\G$-orbits in the Grassmannian of $\QQ$-hyperplanes of $V$ of
signature $(2,n-1)$. Then the corresponding collection of totally
geodesic hypersurfaces $\Hcal_{|\DD}:=\{\DD_H\}_{H\in\Hcal}$ is locally finite on
$\DD$. For a $\QQ$-isotropic line $I$ in $V$, the $H\in\Hcal$
containing $I$ correspond to $\QQ$-hyperplanes in $I^\perp/I\otimes I$
of signature $(1,n-2)$ and so according to example
\ref{example:qcone}-b these decompose $C_{I,+}$ into locally rational
cones. It follows from this and Example \ref{example:qcone}-b that
conditions $(i)$ and $(ii)$ of Definition \ref{def:semitoricsystem} are
satisfied: the space $J_\Hcal$ associated to a $\QQ$-isotropic plane
$J$ is the common intersection of $J^\perp$ and the hyperplanes
$H\in\Hcal$ containing $J$. We denote this admissible partition by
$\Sigma (\Hcal)$.

Notice that for empty $\Hcal$, we recover the Baily-Borel decomposition.
\end{example}

Fix an admissible $\G$-invariant decomposition $\Sigma$ of the conical
locus of $\DD$. We shall see that this determines a compactification of
$X=\G\bs \DD$. We first introduce some notation.

\begin{notation}\label{not:sigma} 
We abbreviate $\pi_{V_\sigma}$ (which we recall,
stands for the projections $V\to V/V_\sigma$ and $\PP
(V)-\PP(V_\sigma)\to \PP (V/V_\sigma)$) by $\pi_\sigma$. We denote the
$\G$-stabilizer of $\sigma$ by $\G_\sigma$ and by $\G^\sigma$ the
subgroup of $\G_\sigma$ that acts trivially on $V/V_\sigma$. The group
of \emph{real} orthogonal transformations of $V$ which leave $\sigma$
pointwise fixed and act as the identity on $V_\sigma/\la \sigma\ra_\CC$
and $V/V_\sigma$ is denoted $N_\sigma$. 
\end{notation}

So $N_\sigma$ is a subgroup of
$N_{I(\sigma)}$, the group of real orthogonal transformations of $V$
which act as the identity on $I(\sigma)$, $I(\sigma)^\perp/I(\sigma)$
and (hence) $V/I(\sigma)^\perp$. 
To be concrete, if $I=I(\sigma)$ is a line, then $N_\sigma$ is the
subgroup of $N_I$ that is the image of $I\otimes V_\sigma/I)(\RR)$
under the exponential map; if it is a plane, then $N_\sigma$ is a
Heisenberg subgroup of $N_I$ that modulo its center is the image of
$(I\otimes V_\sigma/I)(\RR)$ under the exponential map. We have $N_{
\{0\}}=\{ 1\}$, of course.

Observe that $N_\sigma\exp (\sqrt{-1}\sigma)$ is a semisubgroup of
$\Orth (\phi)$ which leaves each fiber of $\DD\to\pi_\sigma\DD$
invariant. (In fact, it is not difficult to verify that every such fiber is an orbit 
of the semigroup that we get if we replace in the definition $\sigma$ by its relative
interior.)

We form the disjoint union 
\[ 
\DD^\Sigma :=\coprod_{\sigma\in \Sigma}\pi_\sigma\DD 
\] 
and define a topology on it as follows. For a subset
$K\subset\DD$ and $\sigma\in\Sigma$ we define its
\emph{$(\G,\sigma)$-saturation} in $\DD^\Sigma$ by 
\[ 
K^\Sigma(\G,\sigma):= \coprod_{\tau\in \Sigma ,\tau\leq \sigma} 
\pi_\tau(\G^\sigma N_\sigma \exp(\sqrt{-1}\sigma)K). 
\] 
These saturated subsets form the
basis of a topology if $K$ runs over the open subsets of $\DD$ and
$\sigma$ over $\Sigma$. We regard $\DD^\Sigma$ as a ringed space in the
usual way: it comes equipped with the sheaf of complex valued, 
continuous, piecewise analytic functions. Similarly for
$(\LL^\times)^\Sigma :=\sqcup_{\sigma\in \Sigma} \pi_\sigma\LL^\times$.
It is clear that in either case $\G$ acts as a group of homeomorphisms.
We put $X^\Sigma:=\G\bs\DD^\Sigma$ (as a ringed space). The image
of $\pi_\sigma\DD$ is the orbit space of the latter by the group
$\G (\sigma):=\G_\sigma/\G^\sigma$; we denote this orbit space by 
$X(\sigma)$. This is an algebraic torus when $\dim I(\sigma)=1$ and 
the total space of a bundle of abelian varieties over a modular
curve when $\dim I(\sigma)=2$.

Notice that we recover the Baily-Borel extensions $\DD^\bb$ and $X^\bb$
as  topological spaces if we take for $\Sigma$ the coarsest admissible 
decomposition. We shall see that as in the Baily-Borel case, $X^\Sigma$ 
is a normal analytic space.

The proof of the following lemma is straightforward.

\begin{lemma} Let $\Sigma'$ be a $\G$-invariant admissible
decomposition of the conical locus that is refined by $\Sigma$. Given
$\sigma\in\Sigma$, let $\sigma'$ be the smallest member of $\Sigma'$
that contains $\sigma$. Then $V_\sigma\subset V_{\sigma'}$, so that we
have a factorization of $\DD\to\pi_{\sigma'}\DD$ over $\pi_\sigma\DD$.
The disjoint union of the maps $\pi_\sigma\DD\to\pi_{\sigma'}\DD$ over
all members $\sigma\in\Sigma$ defines a continuous, $\G$-equivariant
morphism of ringed spaces $\DD^\Sigma\to \DD^{\Sigma'}$. In particular
we have a morphism of ringed spaces $X^\Sigma\to X^{\Sigma '}$.
Likewise for $(\LL^\times)^\Sigma\to (\LL^\times)^{\Sigma'}$. \end{lemma}

This applies in particular to the case when $\Sigma'$ is the
Baily-Borel decomposition so that we always have a map $X^\Sigma\to
X^\bb$.

\begin{theorem} 
The ringed space $X^\Sigma$ is a normal analytic space.
The pull-back of the automorphic $\CC^\times$-bundle on $X^\bb$ to
$X^\Sigma$ can be identified with $(\G\bs \LL^\times)^\Sigma$.
\end{theorem} 
\begin{proof} 
We need to verify this over a neighborhood
of any boundary stratum $X(I^\perp)$ of $X^\bb$. According to the discussion
in \ref{bbext} this then becomes an issue on the level of a Baily-Borel
extension $\G^I\bs Star(\pi_{V^I}\DD)$. The two cases (corresponding to
$\dim I=2$ and $\dim I=1$) are then covered by Section \ref{section:onebc} 
and  Proposition \ref{prop:analyticspace}.
The last assertion is easy. 
\end{proof}

\section{Arrangements on a type IV domain}\label{section:arriv}
We now assume given a
$\G$-invariant arithmetic arrangement $\Hcal$ in $V$ as in Example
\ref{example:qcone}-b. We write 
\[ 
\DD^\circ :=\DD-\cup_{H\in\Hcal}\DD_H \text{  and }\quad X^\circ
:=\G\bs\DD^\circ \] 
for the arrangement complement and its orbit space
respectively. We shall describe a compactification of $X^\circ$ that is
similar to what we did in \cite{I} for the ball quotient case.
We noted in  Example \ref{example:arrdec} that $\Hcal$ defines an
admissible partition $\Sigma (\Hcal)$ of the conical locus. 
Any $H\in\Hcal$ defines a
hypersurface $X_H$ in $X$ which only depends on the image of $H$ in
$\G\bs\Hcal$. This collection hypersurfaces 
is finite. We also regard $X_H$ as a Weil divisor $X_{H}$ on $X^\bb$ 
and we denote its strict transform in $X^{\Sigma (\Hcal)}$ by 
$X_H^{\Sigma (\Hcal)}$.

Analogous to the proof of Lemma (5.1) of \cite{I} one establishes:

\begin{lemma} After passing to a subgroup of $\G$ of finite index, each
hypersurface $X_H$ is without self-intersection in the sense that
its normalization is a homeomorphism that is an isomorphism over $X_H$.
\end{lemma}

From now on we assume that we have made this passage, so that each
$X_H$ is without self-intersection.

A linear subspace $L$ of $V$ contains a vector $z$ with $\phi (z,z)>0$
if and only if $\PP(L)$ meets $\DD$. In that case $\DD_L:=\DD\cap
\PP(L)$ is totally geodesic in $\DD$ and is the bounded symmetric
domain for its group of complex-analytic automorphisms. If moreover $L$
is defined over $\QQ$, then the $\G$-stabilizer acts as an arithmetic
group on $L$ and the closure of $\DD_L$ in $\DD^\bb$ is just the
Baily-Borel extension of $\DD_L$ and hence we denote that closure
$\DD_L^\bb$. It is not hard to see that a similar property holds for
the closure $\DD_L$ in $\DD_L^{\Sigma(\Hcal)}$: is it the extension defined by
the restriction of $\Hcal$ to $L$; we denote that closure by
$\DD_L^{\Sigma (\Hcal )}$.

Let $\PO(\Hcal_{|\DD})$ denote the collection intersections $L$ from
members of $\Hcal$ which have the above property, i.e., which contain a 
vector
$z$ with $\phi (z,z)>0$. For an $L$ of this type, the image $X_L$ of
$\DD_L$ in $X$ depends only on the image of $L$ in $\G\bs
\PO(\Hcal_{|\DD})$. It is a connected component of intersections of
members of the arrangement $\G\bs \Hcal$ (and each such component occurs 
once
and only once). The closure of $X_L$ in $X^\bb$ resp.\ $X^{\Sigma (\Hcal)}$ is the
image of $\DD_L^\bb$ resp.\ $\DD_L^{\Sigma (\Hcal )}$ and denoted $\bar X_L$ resp.\
$X_L^{\Sigma (\Hcal )}$.

\begin{proposition}
The closure $X^{\Sigma (\Hcal)}_H$ of the hypersurface $X_H$ in $X^{\Sigma (\Hcal)}$ supports
an effective Cartier divisor and the collection of these Cartier divisors
is an arrangement on $X^{\Sigma (\Hcal)}$ linearized by $\Lcal$.
If no one-dimensional intersection from $\Hcal$ is positive semidefinite,
then $X^{\Sigma (\Hcal)}\to X^\bb$ is the normalized blowup of the reduced sum of these Cartier divisors.
\end{proposition}
\begin{proof}
That closure of $X_H$ in $X^{\Sigma (\Hcal)}$ is Cartier follows from Lemma's
\ref{lemma:cartier1} and \ref{lemma:cartier2}.
The condition that no one-dimensional intersection from 
$\Hcal$ is positive semidefinite amounts to the nonoccurence
in $X^{\Sigma (\Hcal)}$ of strata of codimension one. So the last assertion
follows form  Corollaries \ref{cor:blowup1} and Corollary\ref{cor:blowup2}.

As for what remains, let $H\in\Hcal$ be defined by the 
$\QQ$-linear form $f:V\to\CC$. We regard $f$ as an object 
on $\DD$ in a standard manner, namely as a section of 
$\LL^{-1}$ with divisor $\DD_H$. This
section is clearly equivariant with respect to the action
of the $G$-stabilizer $G_f$ of $F$. It extends across the union
of the strata of $\DD^{\Sigma (\Hcal )}$ that are of the form 
$\pi_W\DD$, where $W$ is a linear subspace of $H$ (simply because
$f$ factors through $V/W$ in that case). The union of these strata
is open and makes up a neighborhood of the closure 
$\DD_H^{\Sigma (\Hcal )}$ of $\DD_H$ in $\DD^{\Sigma (\Hcal )}$. 
The extension in question is still $G_f$-invariant and has 
$\DD_H^{\Sigma (\Hcal )}$ as `divisor'. This drops and restricts 
to a generating section of the coherent restriction 
of $\Lcal^{-1} (-X_H)$ to $X^{\Sigma (\Hcal)}_H$.
\end{proof}

\begin{proposition}\label{prop:merforms} 
Given $H\in\Hcal$, then the sheaf $\Lcal^{\otimes
k}(kX^{\Sigma (\Hcal)}_H)$ on $X^{\Sigma (\Hcal)}$ is generated by its sections
when $k$ is large enough. 
\end{proposition}
\begin{proof}
Let $v_0\in V(\QQ)$ be 
perpendicular to $H$ and denote by $\Ocal$ the $\G$-orbit of $v_0$. 
We claim that for $k\ge n+2$, the series
 \[ 
F^{(k)}(z):=\sum_{v\in\Ocal} \phi(z,v)^{-k} 
\] 
represents a $\G$-invariant meromorphic section of $\LL^{\otimes k}$.
To see this, let $K\subset \LL^\times$ be compact. Denote by $Y$ the
hypersurface in $V(\RR)$ defined by $\phi (x,x)=\phi (v,v)$.
It has real dimension $n+1$. The $\G$-orbit  $\Ocal$ generates a lattice
in $V(\QQ)$. The number $A(N)$ of points $w$ of this lattice 
for which $\min_{z\in K} |\phi (z,w)|$ lies in the
interval $[t-1,t]$ is for large $N$ proportional to the
$(n+1)$-dimensional volume of the corresponding subset of $Y$. So
$A(t)$ is bounded by a polynomial of degree $n$ in $t$. We find that in
the series  only finitely many terms have a pole on 
$K$ whereas the sum over the absolute values
of the other terms is bounded by a constant plus $\sum_{t\ge 2} A(t)
(t-1)^{-k}$. This converges since $k\ge n+2$. 
We extend $F^{(k)}$ meromorphically over
$(\LL^\times)^{\Sigma (\Hcal )}$ by letting its restriction to $\pi_{V_I}\LL^\times$
be  \[
\pi_{V_I}(z)\in \pi_{V_I}\LL^\times\mapsto
\sum_{v\in\Ocal\cap  I^\perp} \phi(z,v)^{-k}. 
\] 
A minor modification of the arguments used in the proof
of Lemma \ref{lemma:cartier1} and Proposition \ref{prop:semitoricsep} 
shows that this
extension indeed defines  a section of $\Lcal^{\otimes k}(kX_H^{\Sigma (\Hcal )})$.
Notice that $\pi_{V_I}(z_o) \in\pi_{V_I}\LL^\times$
lies in the polar locus of $F^{(k)}$ if and only if 
there exists a $v\in\Ocal\cap I^\perp$ with $\phi (z_o,v)=0$. 
This $v$ is the unique (since we assumed $X_H$ to be without 
selfintersection) and the polar part of $F^{(k)}|\pi_{V_I}\DD$ at 
$\pi_{V_I}(z_o) \in\pi_{V_I}\LL^\times$ 
is the single term $\phi(z,v)^{-k}$. 

It is well-known that for $k$ large enough, $\Lcal^{\otimes k}$ is
generated by its  sections (in fact the same construction---but 
now with $\phi (v_o,v_o)\le 0$---produces enough of these).
The proposition follows.
\end{proof}

Let $\widetilde{X^\circ}\to X^{\Sigma (\Hcal)}$ the blowup of the arrangement defined
by $\Hcal$ in the sense of \cite{I}. This means that we blow up
successively the connected components of intersections of members of
the arrangement, in the order of increasing dimension. We prefer to
obtain this as the $\G$-orbit space of a similar construction on
$\DD^{\Sigma (\Hcal )}$: 
$\widetilde{\DD^\circ}\to \DD^{\Sigma (\Hcal )}$ is the successive blowup of
the submanifolds $\DD_L^{\Sigma (\Hcal )}$.

We denote the exceptional hypersurface associated to 
$L\in\PO(\Hcal_{|\DD})$
by $E(L)\subset\widetilde{\DD^\circ}$. The linearization defines a projection
$E(L)\to\PP (V/L)$. The members of $\Hcal$ that contain $L$ are finite
in number and define an arrangement on $\PP (V/L)$. The arrangement
complement in the latter is denoted $\PP (V/L)^\circ$. The preimage
$E(L)^\circ$ of this arrangement complement in $E(L)$ is naturally a
product: $L\times\PP (V/L)^\circ$ and so we have a projection
$E(L)^\circ\to \PP (V/L)^\circ=\pi_L\DD^\circ$.

The hypersurfaces indexed by a subset of $\PO(\Hcal_{|\DD})$ have a 
nonempty
common intersection if and only if that subset is linearly ordered,
i.e., of the form $L_0\subset L_1\subset\cdots\subset L_r$. The
corresponding intersection $E(L_0,\dots ,L_r)=\cap_i E(L_i)$ is the
closure of a stratum that we will denote by $E(L_0,\dots ,L_r)^\circ$.
Now $\widetilde{\DD^\circ}$ is the disjoint union of 
the arrangement complement $(\DD^{\Sigma (\Hcal)})^\circ$ and the collection $E(L_\pt)^\circ$, where $L_\pt$ runs over the
nonempty linearly ordered subsets of $\PO(\Hcal_{|\DD})$. A quotient set
$\widehat{\DD^\circ}$ of $\widetilde{\DD^\circ}$ is defined by the disjoint
union of $(\DD^{\Sigma (\Hcal)})^\circ$ and the projective arrangement complements $\PP (V/L)^\circ$ with the quotient map
$\widetilde{\DD^\circ}\to\widehat{\DD^\circ}$ being the union of the identity
on $(\DD^{\Sigma(\Hcal)})^\circ$ and the projections $E(L)^\circ\to\PP
(V/L)^\circ$. We give $\widehat{\DD^\circ}$ the quotient topology.
Let us see what happens over a boundary component.

A $\QQ$-isotropic plane $J$ defines the stratum $\pi_{J^\Hcal}\DD$ of
$\DD^{\Sigma (\Hcal )}$, where we recall that $J^\Hcal$ is the intersection of
$J^\perp$ with all the members of $\Hcal$ that contain $J$. The subset
$\Hcal^J$ of $\Hcal$ defined by this last condition defines an
arrangement in $\pi_{J^\Hcal}\DD$. This arrangement is blown up and
then blown down in the manner we described earlier. The arrangement
complement is clearly unaffected by this and so appears in
$\widehat{\DD^\circ}$; it is just $\pi_{J^\Hcal}\DD^\circ$. The total
transform of the arrangement on $\pi_{J^\Hcal}\DD$ is already accounted
for, at least in $\widehat{\DD^\circ}$, since it will map to the union
of the strata $\pi_L\DD^\circ$.

When $I$ is a $\QQ$-isotropic line, then the collection $\Hcal^I$ of
members of $\Hcal$ containing $I$ define an arrangement on
$\DD\cong\pi_I\DD$. The latter yields a semitoric embedding. A stratum
of this torus embedding that does not come from an isotropic plane is
defined by a  $\sigma\in\Sigma$ which meets $C_I$. In that case
$\pi_\sigma\DD=\pi_{V_\sigma}\DD$, where $V_\sigma$ is the common
intersection of $I^\perp$ and the members of $\Hcal$ that contain $I$.
We conclude that 
\[ 
\widehat{\DD^\circ} =\DD^\circ \sqcup
\coprod_{L\in \PO(\Hcal_{|\DD})} \pi_L\DD^\circ \sqcup \coprod_{\sigma\in
\Sigma (\Hcal)}\pi_{V_\sigma}\DD^\circ . 
\] 
In this disjoint union we may have repetitions: this happens for instance
if two distinct members of $\Sigma$ which meet $C_I$ have the same linear span. Since $\DD^\circ$ is open
and dense in $\widehat{\DD^\circ}$, we regard the latter as an
extension of the former. If we pass to the $\G$-orbit space, then 
\[
\widehat{X^\circ}:=\G\bs\widehat{\DD^\circ} 
\] 
is a compactification of
$X^\circ$ whose boundary is decomposed into finitely many strata. Each
boundary stratum has the structure of an arrangement complement, which
can be of projective type, of relative abelian type (over a modular
curve) or of toric type.

As in part I \cite{I}, the $k$-th power of a fractional
ideal $\Ical$ on a variety $Y$ (where $k$ is a positive integer),
denoted $\Ical^{(k)}$, is the image of $\Ical^{\otimes k}$ in the sheaf
of rational functions on $Y$ and hence itself a fractional ideal. We
put 
\[ 
\Lcal (\Hcal) :=\sum_{H\in\G\bs\Hcal}
\Lcal_{X^{\Sigma (\Hcal)}}(X_H^{\Sigma (\Hcal )}). 
\]
Since the arrangement on $X^{\Sigma (\Hcal)}$ is linearized by $\Lcal$, the 
coherent pull-back of $\Lcal (\Hcal)$ to $\widetilde{X^\circ}$ 
is an orbiline bundle which is constant along the fibers of the map 
$\widetilde{X^\circ}\to\widehat{X^\circ}$. 

\begin{theorem}\label{thm:contract} 
The coherent pull-back of $\Lcal
(\Hcal)$ to $\widetilde{X^\circ}$ is a semi-ample line bundle which
defines the contraction $\widetilde{X^\circ}\to\widehat{X^\circ}$ (by which
we mean that a suitable power of this bundle defines a morphism onto a
projective variety which as a topological quotient of $\widetilde{X^\circ}$
is just $\widehat{X^\circ}$). 
\end{theorem}
\begin{proof}
That the coherent pull-back of $\Lcal (\Hcal)$ to $\widetilde{X^\circ}$ 
is an line bundle is a general property of the blowup of an
arrangement. Its $k$th power is generated as a line bundle
by the subbundles $\Lcal^{\otimes k}(X_H^{\Sigma (\Hcal )})$. Since these subbundles
are generated by their sections for large $k$, the 
same is true for large powers of $\Lcal (\Hcal)$.
A section of a power of $\Lcal (\Hcal)$ is constant 
along the fibers of the map $\widetilde{X^\circ}\to\widehat{X^\circ}$.
Since $\Lcal$ is ample, they separate these fibers.
The theorem follows.
\end{proof}

\begin{corollary}\label{maincor:iv} 
Suppose that no one-dimensional intersection of
members of the arithmetic arrangement $\Hcal$ in $V$ 
is positive semidefinite. Then the algebra of automorphic forms
\[ 
\oplus_{k\in\ZZ} H^0(\DD^\circ,\Ocal (\LL^k))^\G 
\] 
(where $\LL$ is the natural automorphic bundle over $\DD$ and $\DD^\circ$
is the arrangement complement) is finitely generated with positive degree
generators and its proj is the compactification $\widehat{X^\circ}$ of
$X^\circ$. Moreover, the boundary $\widehat{X^\circ}-X^\circ$ is the
strict transform of the boundary $X^{\Sigma (\Hcal)}-X$. 
\end{corollary}
\begin{proof}
The assumption implies that $(X^{\Sigma (\Hcal)})^\circ-X^\circ$
is dense in $\widehat{X^\circ}-X^\circ$. The last assertion then follows 
and for the other assertions it suffices to note that 
$X^{\Sigma (\Hcal)}-X$ is of codimension $\ge 2$ in $X^{\Sigma (\Hcal)}$.  
\end{proof}

In applications pertaining to $K3$-surfaces, we often
get the finite generation of the algebra of automorphic forms 
via the following theorem. The previous corollary can then be used 
to interpret its proj.

\begin{theorem}\label{thm:gitcon} 
Let $Y$ be a normal projective
variety with ample line bundle $\eta$, $G$ a reductive group  acting
on the pair $(Y,\eta)$ and $U\subset Y$ a $G$-invariant open-dense
subset on which $G$ acts properly. 
Suppose that any $G$-invariant section of a tensor power of 
$\eta$ over $U$
extends across $Y$ (a condition that is certainly satisfied if $Y-U$ is of 
codimension
$>1$ in $Y$). Let $\Hcal$ be an arithmetic arrangement on a type IV 
domain $\DD$ relative to an arithmetic group $\G$. 

Then an identification (in the analytic category) of the line 
bundle $G\bs (U, \eta |U)$ with the pair $(X^\circ,\Lcal |X^\circ)$
associated to this  arrangement gives rise to an isomorphism of 
$\CC$-algebra's
\[ 
\oplus_{k=0}^\infty H^0(Y,\eta^{\otimes k})^G\cong
\oplus_{k\in\ZZ} H^0(\DD^\circ,\Ocal (\LL)^{\otimes k})^\G. 
\] 
In particular, the algebra of meromorphic automorphic forms is finitely
generated with positive degree generators, $U$ consists of stable 
$G$-orbits, and 
$G\bss Y^\ss\cong \widehat{X^{\circ}}$.
\end{theorem} 
\begin{proof} 
To say that we have an isomorphism between 
$H^0(G\bs U,G\bs (\eta^{\otimes k}|U))$ and 
$H^0(X^\circ,\Lcal^{\otimes k})$ is saying that we have an isomorphism
between $H^0(U,\eta^{\otimes k})^G$ and
$H^0(\DD^\circ,\Ocal(\LL)^{\otimes k})^\G$.
Since $H^0(Y,\eta^{\otimes k})^G\to H^0(U,\eta^{\otimes k})^G$ is 
an isomorphism by assumption, the assertion follows.
\end{proof}

As in the ball quotient case, the identification demanded by the
theorem will often come from a period mapping. In the following sections
we shall discuss a few examples related to $K3$-surfaces.

Borcherds  described in \cite{borcherds}  a technique that produces
in a systematic fashion arithmetic arrangements that are definable by an 
automorphic form. Here we note a necessary condition for this to be the case: 
the divisor on $X$ defined by such an arrangement 
must then extend across the Baily-Borel compactification $X^\bb$ as a Cartier 
divisor.  Proposition \ref{prop:weylv} gives us immediately:

\begin{theorem}
Let a divisor on $X$ be defined by the $\G$-invariant arrangement 
$\Hcal$ and the $\G$-invariant function $m:\Hcal\to\ZZ$. Then the divisor
has the Cartier property at the boundary if and only
if  for every $\QQ$-isotropic line $I$ there exists a quasi-character 
$e^{\rho_I}: \PP (V)-\PP (I^\perp)\to\CC$  satisfying the functional 
equation of 
Proposition \ref{prop:weylv} relative to $\G_I$ and the restriction of
$m$ to the arrangement defined by the collection of $H\in\Hcal$ containing $I$.
\end{theorem}

It is clear that  this property needs only be tested 
for a system of representatives of the $\G$-orbits in the collection of  
$\QQ$-isotropic lines.  Still, it is quite a strong condition. 
If it is satisfied and if the Picard group of $X^\bb$ happens to have rank one, 
then it follows that a multiple of the corresponding divisor is defined by an 
automorphic form. It would be nice to have here an analogue of Proposition 
\ref{prop:weylv}.

\section{Applications to period mappings of polarized $K3$-surfaces}

\subsection{Semiample linear systems on $K3$-surfaces}
We return to the moduli spaces of $K3$-surfaces. Suppose $S$ is a
$K3$-surface equipped with a primitive semi-ample class $h\in\pic (S)$
with $h\cdot h=2g-2$. The linear system defined by $h$ consist of genus
$g$ curves and is of dimension $g$. According to A.~Mayer \cite{mayer},
we are then in one of the following situations: 
\begin{itemize}
\item[---] The \emph{unigonal} case: the linear system maps to a
rational curve. This happens precisely when there is an elliptic curve
on $S$ on which $h$ has degree $1$, a condition equivalent to the
existence of an isotropic  class $f\in \pic(S)$ with $h\cdot f=1$.
The moving part of the linear system gives $S$ the structure
of an elliptic surface with section.  
\item[---] The \emph{digonal} case:
the linear system defines a morphism of degree $2$ on its image
(a rational surface). This is always so when $g=2$ 
(unless we are in the unigonal case, of course): 
then the linear system realizes $S$ as a
double cover of $\PP^2$ ramified along a sextic curve with simple
singularities. When $g>2$, this happens precisely when there is an
elliptic curve on $S$ on which $h$ has degree $2$, a condition
equivalent to the existence of a primitive isotropic class $f\in
\pic(S)$ with $h\cdot f=2$, again assuming we are not in the unigonal
case. The moving part of the linear system gives $S$ the structure
of an elliptic surface, in general without a section. 
\item[---] The \emph{nonhyperelliptic} case: the linear system maps
$S$ birationally onto its image. 
\end{itemize} 
This means that the question of whether a given pair $(S,h)$ is 
unigonal or digonal can be
read off from its periods: Adhering to the notation in the final
subsection of Section \ref{section:onebc},
consider the set $E_g'$ resp.\  $E_g''$ of primitive isotropic vectors
$f$ in the $K3$-lattice $\Lambda$ with $f\cdot h_g=1$ resp.\ $f\cdot h_g=2$. It is
well-known that each of these sets is a single orbit of $\G_g$. 
Representative elements are $f_3$ and $2f_3+f_2$ respectively.
The span of $h_g=e_3+(g-1)f_3$ and $f_3$ is a hyperbolic summand $U$,
whereas the span of $h_g$ and $2f_3+f_2$ is
isomorphic to $I(2)\perp I(-2)$ or $U(2)$,
depending on whether $g$ is even or odd. These are
even rank two lattices of hyperbolic signature so that their
orthogonal complements have signature $(2,18)$. We notice in passing
that in each case the discriminant group of the rank two lattice 
is $2$-torsion, which implies that this lattice 
is the fixed point set of an involution in $\G_g$. Such an 
involution acts on $\DD_g$ as reflection.
So the collection $\Hcal'_g$ resp.\ $\Hcal''_g$ of
hyperplanes in $\Lambda_g\otimes\CC$ perpendicular to the span 
of $h_g$ and a vector in
$E_g'$ resp.\  $E_g''$ make up an `reflection arrangement' on $\DD_g$. 
Each of these in turn determines a divisor $D_g'$ resp.\
$D_g''$ in $X_g$. Via the period mapping, the generic point of $D'_g$ resp.\
$D''_{g\ge 3}$ corresponds to the property of the $K3$-surface being
unigonal resp.\ digonal. We therefore put
\[
E_g:=
\begin{cases} 
E'_2\text{ if $g=2$},\\
E'_g\cup E''_g \text{ if $g\ge 3$.}
\end{cases}
\]
We denote by $\Hcal_g$ the corresponding collection of hyperplanes
and by $X_g^\circ$ the associated arrangement complement.
We will see that for small $g$, $\widehat{X_g^\circ}$ has the interpretation
of a  GIT compactification. In order to understand the strata of 
$X_g^{\Sigma (\Hcal_g)}$ we need the following lemma.

\begin{lemma}\label{lemma:k3arr}
Let $M\subset\Lambda$ be a nondegenerate sublattice of signature 
$(1,r)$ with $r\ge 2$ spanned by $h_g$ and elements of $E_g$. 
Then $g\in\{ 3,4\}$, $M\cap E_g=M\cap E''_g$ and 
and $v_1\cdot v_2=1$ for every pair $v_1,v_2\in M\cap E''_g$
which spans with $h_g$ a rank $3$ lattice.
In case $g=4$, $M$ is isometric to $U\perp I(-2)$ (and so $r=2$)  and 
$M\cap E''_g$ consists of three elements whose sum is $h_4$.
\end{lemma}
\begin{proof}
Suppose $v_1,v_2\in M\cap E_g$ are such that they form with $h_g$
a linearly independent set. Then the
quadratic form on the span of $h_g,v_1,v_2$ has in this basis
the expression $(2g-2)x^2+2a_1xy_1 +2a_2xy_2 +2\lambda y_1y_2$, 
with $a_i\in\{1,2\}$, depending on whether $v_i$ is in $E_g'$ or $E_g''$.
The determinant  of this  quadratic form is
$2\lambda (a_1a_2-(g-1)\lambda)$. We want it to be negative, which
means that $0<\lambda < a_1a_2(g-1)^{-1}$. As $\lambda$
is an integer, this can only be if $g-1\in\{ 2,3\}$, $a_1=a_2=2$
and $\lambda=1$ (recall that when $g=2$, we have $a_i=1$, 
by assumption). This proves the first part of the lemma.

If $v_1,v_2,v_3$ are three distinct elements of $M\cap E''_4$, then 
the lattice spanned by the elements $v_1,v_2,v_3$
is isomorphic with $U\perp I(-2)$ (take as basis $(v_1, v_2, v_3-v_1-v_2)$
and so has hyperbolic signature.
Now notice that $h_4-v_1-v_2-v_3$ is perpendicular to the sublattice 
$N\subset M$ spanned by $v_1, v_2, v_3, h_4$. Since $M$ has hyperbolic
signature, this can only be if $h_4=v_1+v_2+v_3$. 
It also follows that $M\cap E''_4$ cannot have more than three  elements. 
\end{proof}

\begin{remark}
Let us briefly comment on the remaining (and quite interesting) case $g=3$.   
If $v\in  M\cap E''_3$, then its image under the orthogonal projection 
$\Lambda\to \Lambda_3\otimes\QQ$ is $v-\tfrac{1}{2}h$, whose 
double $2v-h_3$ is a $(-4)$-vector. Notice that $h_3-v\in  M\cap E''_3$ also
and that its image in $\Lambda_3\otimes\QQ$ is the antipode of
the orthogonal projection of $v$. If $v'\in M\cap E''_3$ is distinct from
$v$ and $h_3-v$, then the orthogonal projections $v-\tfrac{1}{2}h$ and
$v'-\tfrac{1}{2}h_3$ are perpendicular. This shows that $M\cap E_3''$
has exactly $2r$ elements $v_1, v_2,\dots ,v_r, h_3-v_1,\dots ,h_3-v_r$
whose orthogonal projection in $\Lambda_3\otimes\QQ$
is the union of an orthonormal set and its antipode. 
Let us denote by $M_r$ the abstract lattice defined by such a system 
$(h_3, v_1,\dots ,v_r)$  and remember that  $M_r$ has the distinguished
vector $h_3$.  Suppose $r\ge 3$. Then the elements $v_1,v_2$ span a 
copy of $U$ in $M_r$ and the orthogonal complement
of this copy in $M_r$ has the basis 
$(h_3-v_1-v_2-v_3, v_1+v_2-v_3,v_3-v_4, v_4-v_5,\dots ,v_r-v_{r-1}$,
which is a root basis of type $D_{r-1}$. So $M_r\cong U\perp D_{r-1}(-1)$.  
In order that $M_r$ embeds in $\Lambda$, we must have $r\le 19$, for
reasons of signature. The value $r=19$ is attained, 
for one can easily show that 
$U\perp D_{18}(-1)\perp I(2)\perp I(2)$ has an even unimodular overlattice 
(of signature $(3,19$) and as is well-known, such a lattice is isomorphic to $\Lambda$. 
\end{remark}

\begin{lemma}\label{lemma:uniqueemb}
The group $\G_g$ acts transitively on $E'_g$ for $g\ge 2$) 
and on $E''_g$ for $g\ge 3$. If $g\ge 4$ and $r$ is a
positive integer, then $\G_g$ is also transitive on the collection of 
nondegenerate sublattices of  $\Lambda$ of signature $(1,r)$ spanned 
by $h_g$ and elements of $E''_g$, provided that this collection is nonempty. 
\end{lemma}
\begin{proof}
We will first prove this lemma for equivalence with respect to the stabilizer 
$\Orth (\Lambda)_{h_g}$ of $h_g$ in the orthogonal  group of $\Lambda$. 
For this we need the following well-known \emph{transitivity property}: given 
an even nondegenerate lattice $M$ of rank $\rho$ and an even unimodular 
lattice 
$N$ which contains a copy of $U^{\perp (\rho+1)}$, then the primitive 
embeddings 
of $M$ in $N$ are all in the same orbit of the orthogonal group of $N$. 

For all $g\ge 2$ resp.\ $g\ge 3$,  $h_g$ and an element of $E'_g$ resp.\ 
$E''_g$ span a primitive sublattice of  $\Lambda$
and so the above transitivity property implies that $E'_g$ ($g\ge 2$) and 
$E''_g$ ($g\ge 3$) is an orbit of $\Orth (\Lambda)_{h_g}$. 
This already covers all cases $g\ge 5$. 
For $g=4$, the remaining  case is when the sublattice in question is isometric 
to $U\perp I(-2)$.
Since $U$ embeds uniquely up to orthogonal transformation in $\Lambda$, the 
issue is equivalent to the uniqueness of the embedding of $I(-2)$ in the even 
unimodular lattice $\Lambda':=U^{\perp 2}\perp E_8(-1)^{\perp 2}$. 
This is also covered by the quoted transitivity property.
To finish the proof we observe that in all cases we can find an element in
$\Orth (\Lambda)_{h_g}$ which stabilizes an element of $E'_g$ resp.\ a 
sublattice as in the lemma and sends $\DD_g$ to its complex conjugate. 
\end{proof}

This lemma shows in particular that $D'_g$ and $D''_{g\ge 3}$ are irreducible.

\begin{remark}\label{remark:g=3}
The second clause of the lemma above is not true as stated when $g=3$, since  for $r=9$ and $r=17$ there exist imprimitive embeddings of
$M_r\cong U\perp D_{r-1}$ in $\Lambda$ for which the image of $h_3$ is still primitive.
However it can be shown that Lemma \ref{lemma:uniqueemb} still holds for 
primitive sublattices. The imprimitive 
sublattices---which exist only for $r=9$ and $r=17$---also form a single 
$\G_3$-orbit.
\end{remark}

From Lemma \ref{lemma:k3arr} and  Corollary \ref{maincor:iv} we deduce:

\begin{corollary}
For all $g\ge 2$,  $g\not= 3$, the algebra of meromorphic automorphic forms
\[
\oplus_{k\in\ZZ} H^0(\DD_g^\circ, \Ocal (\LL^k))^{\G_g}
\]
is finitely generated with positive degree generators.
The boundary of the semitoric compactification 
$X_g^{\Sigma (\Hcal_g)}$ of $X_g$ and of the compactification 
$\widehat{X_g^\circ}$ of $X_g^\circ$ is of dimension $2$, except when 
$g=4$, where this dimension is $3$.
The closures of $D_g'$ and $D_g''$ are disjoint in $X_g^{\Sigma (\Hcal )}$ and
$D_g'$ is without selfintersection. The same is true of $D_g''$ when $g>4$.
\end{corollary}

We shall see that the first part of the statement of this corollary also holds for
$g=3$. 

Saint-Donat \cite{stdonat} made a detailed study of the projective models
of $K3$-surfaces and it is interesting to observe that the loci $D_g'$, $D_g''$ and their 
selfintersections parametrize semipolarized $K3$-surfaces that 
appear in his classification.  This correspondence is as follows
(we omit the proofs, which are not difficult):
\begin{enumerate}\item[]
\begin{enumerate}
\item[$D_g'$:] The $K3$-surfaces endowed with an elliptic fibration with a 
section; 
the polarization is $(g-1)$ times the class of the fibration plus the class
of a section.  The morphism to $\PP^g$ defined by the polarization 
is the composite of the fibration and a degree $g$ embedding of the base 
curve in $\PP^g$. These families become pairwise isomorphic if we forget their 
polarization.
\item[$D_{g\ge 3}''$:]  The $K3$-surfaces $S$ of genus $g\ge 3$
endowed with an elliptic fibration and a bisection. The bisection is
a smooth curve of genus $i=0$ or $i=1$ depending on whether $g$ is even or odd.
The polarization is  $[\frac{1}{2}g]= \frac{1}{2}(g-i)$ times the 
class of the fibration plus the class of the bisection.  
The morphism $S\to \PP^g$ defined by the polarization 
factors as a degree two cover $S\to\Sigma_i$ (here $\Sigma_i$ is the 
$\PP^1$-bundle over $\PP^1$ with a section of selfintersection $-i$) followed
by the embedding of the latter in $\PP^g$ by the linear system defined by 
a section of selfintersection $-i$ plus $[\frac{1}{2}g]$ times a fiber.
\item[$(D_4'')^{(2)}$:]  The $K3$-surfaces $S$ of genus $4$
endowed with an elliptic fibration, a section $C_0$ and a smooth rational 
curve $C_1$ in a fiber met by the section $C_0$, the  polarization being 
$3$ times the class of the fibration plus the class of $2C_0+C_1$.   
The morphism $S\to \PP^4$ defined by the polarization has image in $\PP^4$
the cone over a smooth rational cubic curve in $\PP^3$; the map to the image
is of degree two and ramifies over the vertex.
\item[$(D_3'')^{(r)}$:] The $K3$-surfaces $S$ of genus $3$
endowed with an elliptic fibration,  a section $C_0$ and
smooth rational curves $C_1,\dots ,C_{r-1}$ such that
the intersection diagram of $(\text{(fiber)}+C_0, C_1,\dots, C_r)$
is the affine Dynkin diagram type $\widehat B_r$: for $r=2$, $C_1$ is a section
disjoint from $C_0$, for $r=3$, $C_1$ and $C_2$ lie in distinct fibers
and both meet $C_0$, for $4\le r\le 19$ $C_1,\dots ,C_{r-1}$ lie
in single fiber and make up with $C_0$ a $D_r$-configuration. The 
polarization is given by the positive generator  of the radical
of the $\widehat B_r$-diagram: 
$2(\text{(fiber)}+C_0+\cdots +C_{r-3})+C_{r-2}+C_{r-1}$.
The morphism $S\to \PP^3$ defined by the polarization has image
a quadric cone; the map to the image is of degree two and does not 
ramify over the vertex. (For $r\ge 4$, the curves $C_1,\dots ,C_{r-1}$ 
lie on a Kodaira fiber of the elliptic fibration.) 
\end{enumerate}
\end{enumerate}

\subsection{$K3$-surfaces of small genus} 
A $K3$-surface $X$ of genus $2$ for which the linear system is 
of digonal type is realized by the linear system as a double
cover of a projective plane ramifying along a sextic curve with
only simple singularities. A $K3$-surface of genus $3$, $4$ or $5$
for which the linear system is nonhyperelliptic
is realized as a quartic hypersurface in $\PP^3$,
resp.\ a complete intersection of bidegree $(2,3)$ in  $\PP^4$, 
resp.\ a complete intersection of three quadrics in  $\PP^5$.

We fix a vector space $W_g$ of dimension $g+1$ and a generator 
$\alpha_g$ of $\wedge^{g+1}W_g$, which we shall view as a translation
invariant $(g+1)$-form on $W_g^*$. For $g=2,3,4,5$ we define a
projective variety $Y_g$ with very ample class $\eta_g$ on which
the group $\SL (W_g)$ naturally acts through its simple quotient
$\PSL(W_g)$: 
\begin{enumerate}\item[]
\begin{enumerate}
\item[$g=2$:] $Y_2$ is the projective space $\PP (\sym^6 W_2)$ of all
sextics in  $\PP (W_2^*)$ and $\eta_2=\Ocal_{\PP (\sym^6 W_2)}(1)$,
\item[$g=3$:] $Y_3$ is the projective space $\PP (\sym^4 W_3)$ of all
quartics in  $\PP (W_3^*)$ and $\eta_3=\Ocal_{\PP (\sym^4 W_3)}(1)$, 
\item[$g=4$:] $Y_4=\PP(\EE)$, where $\EE\to \PP (\sym^2 W_4)$ is the
vector bundle whose fiber over  $[Q]$ is the cokernel of  $W_4\to \sym^3 
W_4$, 
$w\mapsto Q.w$ (this makes $Y_4$ the space of all complete
intersections of bidegree $(2,3)$ in $\PP (W_4^*)$), and
$\eta_4=\Ocal_{\PP(\EE)}(1)\otimes\Ocal_{\PP (\sym^2 W_4)}(1)$, 
\item[$g=5$:] $Y_5$ is the Grassmannian $G_3\sym^2(W_5)$ (in other words,
the space of all complete intersections
of tridegree $(2,2,2)$ in $\PP (W_5^*)$) and $\eta_5$
the line bundle defining  the Pl\"ucker embedding.
\end{enumerate}
\end{enumerate}
We denote by $U_g\subset Y_g$ the subset of complete intersections 
for which the only singular points are simple singularities.

The following theorem implies among other things that the algebra
of meromorphic automorphic forms 
$\oplus_{k\in\ZZ}H^0(\DD^\circ_g,\Ocal(\LL^k))^{\G_g}$
is generated with positive degree generators also when $g=3$.

\begin{theorem}\label{thm:sextic} 
Let $g\in \{ 2, 3,4,5\}$. Then the period mapping defines an isomorphism
$\PSL (W_g)\bs U_g\to X_g^\circ$. This isomorphism lifts naturally to 
an isomorphism between the line bundle $\PSL (W_g)\bs (U_g,\eta_g|U_g)$
and  the line bundle $(X_g,\Lcal)$ when $g\not=2$ and the quotient of $(X_2,\Lcal)$
by the center $\{\pm1\}$ of $\G_2$ in  case $g=2$.
This, in turn, determines an isomorphism of normal projective 
varieties $\PSL (W_g)\bss Y_g^\ss\cong\widehat{X_g^\circ}$ and an 
isomorphism of graded $\CC$-algebra's 
\[
\oplus_{k=0}^\infty H^0(Y_g,\eta^{\otimes k}_g)^{\PSL (W_g)}\cong 
\oplus_{k\in\ZZ}H^0(\DD^\circ_g,\Ocal(\LL^k))^{\G_g}.
\]
\end{theorem}
\begin{proof}
We begin with showing that given a $y\in U_g$, then every nonzero linear
form $F$ on the line $\eta_g(y)$ defines a $K3$-surface 
$S(F)$ together with a generating section $\alpha (F)$ of its dualizing
sheaf in such a manner that the dependence on $F$ is of degree $-1$:
$S(tF)$ is canonically isomorphic to $S(F)$ and under this isomorphism 
$\alpha(tF)$ corresponds to $t^{-1}\alpha (F)$. 
In the case $g=2$, which we discuss first, this has to be qualified
since the objects of interest all come with a nontrivial involution.

So assume that $g=2$ and let us write $\PP$ for $\PP (W_2^*)$.
Let $\EE$ denote the total
space of the line bundle over $\PP$ whose sheaf of sections
is $\Ocal_\PP(-3)$. Its dualizing sheaf $\omega_\EE$ is naturally
identified with the pull-back of $\omega_\PP(3)$. 
The choice of generator for $\wedge^3W_2$ defines a nowhere zero
section of $\omega_\PP(3)$ and hence a nowhere zero
section of $\omega_E$. Denote this section by $\alpha$.

Now let $y\in U_2$ define the sextic curve $C(y)$ in $\PP(W_2^*)$ and let 
$F\in\sym^6 W_2$ be a defining equation for this
curve (equivalently, a nonzero linear form on $\eta_2(y)$). 
We regard $F$ as a regular function on $\EE$ which is homogeneous
of degree $2$. Then the zero set of  $F-1$, closed off in the projective
completion of $E$ is a double cover $S(F)\to\PP$ ramified along
$C(y)$ and the residue of  $(F-1)^{-1}\alpha$ on $S(F)$ is a nowhere
zero section $\alpha (F)$ of its dualizing sheaf $\omega_{S(F)}$. 
This construction is homogeneous of degree $-\frac{1}{2}$ in the following sense:
multiplication by $t$ in $\EE$ sends $S(F)$ to $S(t^{-2}F)$ and under this
isomorphism $\alpha (t^{-2}F)$ corresponds to $t\alpha (F)$. This is also 
reflected by the fact that the involution sends $\alpha (F)$  to $-\alpha (F)$. 
So a complete invariant for $\pm (S(F),\alpha (F))$ is $\alpha (F)^2$,
viewed as a section of either $\omega_{S(F)}^{\otimes 2}$
or $\omega_{\PP}^{\otimes 2}(C(y))$. We then obtain an isomorphism between 
the lines $\eta (y)$ and  $H^0(\omega_{\PP}(C(y)))^{\otimes 2}$.

For $g=3,4,5$ we proceed in a similar
fashion, although the situation is simpler. For instance, a nonzero
linear form on $\eta_5(y)$ is represented by a nonzero decomposable 
element 
$F_0\wedge F_1\wedge F_2$ of $\wedge^3 (\sym^2 W_5)$ whose factors
define $S(y)$ in $\PP (W_5^*)$. The iterated residue of $\alpha_5$ 
relative to the ordered triple $(F_0, F_1, F_2)$ defines a generator of
the dualizing sheaf of the common zero set of this triple in $W_5^*$.
This $3$-form is  $\CC^\times$-invariant, and by taking the residue at 
infinity, we get a generating section $\alpha (F)$ of $\omega_{S(y)}$.

Let $\widehat S(y)$ be the minimal resolution of $S(y)$.
We recall that $H^2(\widehat S(y),\ZZ)$ is a free abelian group and equipped
with its intersection pairing isomorphic to the $K3$-lattice $\Lambda$. The
isomorphism can be chosen in such a manner that the pull-back to $\widehat
S(y)$ of the hyperplane class of $\PP(W^*_g)$ is mapped to $h_g$. Under
such an isomorphism, $H^0(\omega_{\widehat S(y)})$ is mapped to $\LL_g$ or
its complex conjugate. We can (and will) choose the isomorphism such that
the former occurs. It is clear from the definition that $\G_g$
permutes the collection of such isomorphisms simply transitively, so
we get a period morphism on the level of line bundles
(with the usual qualification in case $g=2$): 
\[ 
\PSL (W_g)\bs (U_g,\eta_g|U_g)\to (X_g,\Lcal). 
\] 
The theory of period mappings for $K3$-surfaces tells us that
this morphism is an open embedding whose image 
on the base is just $X_g^\circ$. If we combine this
with the observation that the deleted part
$Y_g-U_g$ is of  codimension $\ge 2$, all assertions follow.
\end{proof} 

\begin{problem}
We have not been able to generalize this theorem to higher values of $g$,
although work of Mukai \cite{mukai} suggests that this could be done at 
least for $g\le 10$:
he observed that in the range $6\le g\le 10$, a general primitively 
semipolarized $K3$ of genus $g$ is still a complete
intersection on a homogeneous projective variety (in fact a linear section
for $7\le g\le 10$). If this produces all
nonhyperelliptic $K3$-surfaces with rational double points in $\PP^g$ that are
primitively polarized by the hyperplane class for these values of $g$,
then Theorem \ref{thm:sextic} extends to that range. 
 \end{problem}

\begin{problem}
The singular objects 
parametrized by $Y_g$  define for $g\le 5$ a hypersurface in $Y_g$. 
This hypersurface is 
defined by a section of a (computable) power of $\eta_g$. Via Theorem 
\ref{thm:sextic} we deduce that  there must exist a meromorphic 
$\G_g$-automorphic form whose divisor is twice the arrangement defined 
by the $(-2)$-vectors in $\Lambda_g$
plus a linear combination of the arrangements defined by $E_g'$ ($g\ge 2$)
and $E_g''$ ($g\ge 3$). We expect this automorphic form to have a product
expansion. We wonder whether this is true for all $g$. Again, Mukai's 
work suggests this to be so for $g\le 10$. A theorem due to
Borcherds, Katzarkov, Pantev and Shepherd-Barron  \cite{bkps}
implies that there is an automorphic form whose zero set is an arrangement
which contains the locus in question.
\end{problem}

\section{Application to moduli of general Enriques surfaces}
In this section we briefly explain how the work of H.~Sterk 
fits in this setting (but his results go well beyond what is
discussed here; see  his two part paper \cite{sterk}).
Let us recall that an Enriques surface is a surface which 
admits a $K3$ surface as an unramified covering of degree two.
Since a $K3$-surface is simply connected, an Enriques
surface $\bar S$ is essentially the same thing as a $K3$-surface $S$
equipped with a fixed point free involution. All of the
cohomology  of an Enriques surface is algebraic and so the same is true
for the cohomology of $S$ invariant under the involution.

Every Enriques surface $\bar S$ admits a \emph{semipolarization of 
degree $2$}, that is, a class $\bar h\in\pic (\bar S)$ 
intersecting every effective class nonnegatively and with 
$\bar h\cdot \bar h=2$. 
Let $\bar h\in\pic (\bar S)$ be such a class. Then its preimage $h$ in 
$S$ defines a digonal semipolarization of degree $4$ on $S$: 
its linear system  defines a morphism of degree two from $\widehat S$ onto 
a nonsingular quadric or a quadric cone in projective three space. 
Its discriminant curve is the intersection of this quadric with 
a quartic surface. Following Horikawa we call $(\bar S,\bar h)$ 
\emph{general}  or \emph{special} according to whether 
the quadric is smooth or not. We shall concentrate here on the general case. 
Let us first do the discussion for general digonal $K3$-surfaces of genus
$3$.

\begin{none}
We start out with two complex vector spaces $W$, $W'$ of
dimension two and abbreviate $\PP:=\PP (W^*)$ and $\PP':=\PP(W'{}^*)$. 
Let $\EE\to \PP\times\PP'$ be the line bundle whose sheaf of sections is 
the exterior tensor product $\Ocal_\PP(-2)\boxtimes\Ocal_{\PP'}(-2)$. 
A straightforward argument shows that the choice of a generator of
$\wedge^2W\otimes\wedge^2W'$ determines a generating section $\alpha$ of
its dualizing sheaf $\omega_\EE$.

The divisors of bidegree $(4,4)$ on $\PP\times\PP'$ are parametrized
by the projective space $Y:=\PP(\sym^4 W\otimes\sym^4 W')$.
Take $\eta:=\Ocal_Y(1)$ and 
let $U\subset Y$ be the locus parametrizing divisors 
with only simple singularities. 
Choose $y\in U$ and let  $F\in \sym^4 W\otimes\sym^4 W'$ be a nonzero
element lying over $y$. Then $F$ defines a regular function on $E$
which is homogeneous of degree $2$. This gives rise to a double
cover $S(F)\to \PP\times\PP'$  and a generating section $\alpha (F)$
of its dualizing sheaf. Its minimal resolution is naturally a $K3$-surface 
endowed with a genus $3$ semipolarization of digonal type. 
As in the case for $K3$-surfaces of genus $2$, we get a natural 
isomorphism between
the fiber  $\eta (y)$ and the set of isomorphism classes of pairs 
$(S(F),\alpha (F))$. The subgroup $G$ of $\PSL (W\oplus W')$ 
which preserves the direct sum decomposition (allowing the 
summands to interchange) acts on $\sym^4 W\otimes\sym^4 W'$
and hence on $(Y,\Ocal_Y(1))$ and the $G$-orbit space of
$(U,\Ocal_U(1))$ can be understood as an open part of the moduli space 
of digonal $K3$-surfaces of genus $3$.

Next assume we are given involutions $i$ in $W$ resp.\ $i'$ in $W'$ 
with one-dimensional eigenspaces. Notice that the involution $(i, i')$
has four fixed points in $\PP\times \PP'$.
The involution $(i, i')$ also acts naturally in $\EE$ (with the four
fibers as fixed point set). We put
\[
Y_E:=\PP((\sym^4 W\otimes\sym^4 W')^{i\otimes i'})
\text{ and } \eta_E:=\Ocal_{Y_E}(1),
\]
and denote by $U_E\subset Y_E$ the open subset parametrizing divisors 
not passing through a fixed point and with simple singularities only.
For a nonzero $F\in \sym^4 W\otimes\sym^4 W')^{i\otimes i'}$ over $U_E$,
the involution $-(i,i')$ restricts to a fixed point free involution 
on the minimal resolution $\widehat S(F)$ of $S(F)$ and so its orbit variety 
is an Enriques surface. 
This Enriques surface comes naturally with an a semipolarization of 
degree two and is general. 
If $G_E$ is the $G$-centralizer of the involution $(i, i')$, then
$G\bs U_E$ can be understood as the coarse moduli space
of semipolarized Enriques surfaces of general type.
\end{none}

\begin{none}
It is well-known that the 
action of a fixed point free involution of a $K3$-surface $S$
induces in $H^2(S;\ZZ)$
an involution equivalent to the involution $\iota$ of $\Lambda$ defined as
follows: if we write $\Lambda_E$ for the \emph{Enriques lattice} 
$E_8(-1)\perp U$ and identify $\Lambda$ with
$\Lambda_E\perp \Lambda_E\perp U$, then $\iota(a,b,c)=(b,a,-c)$. Notice
that the eigenlattices $\Lambda^+$ resp.\ $\Lambda^-$ are isometric to
$\Lambda_E(2)$ resp.\ $\Lambda_E(2)\perp U$. The name of $\Lambda_E$ 
is explained by the fact that an equivariant  
isometry $H^2(S;\ZZ)\cong \Lambda$ identifies 
$H^2(\bar S,\ZZ)$ modulo its torsion  with the lattice
$\Lambda^+ (\tfrac{1}{2})\cong \Lambda_E$. 
Observe that the signature of $\Lambda^-$ is $(2,10)$. 

We put $h^+:=(e+f,e+f,0)\in\Lambda^+$. It is not difficult to verify
that there are exactly two isotropic vectors in $\Lambda^+$ with inner 
product $2$ with $h^+$, namely $e^+:=(e,e,0)$ and $f^+:=(f,f,0)$. If
$h\in\pic (S)$  is the pull-back of an semipolarization of degree 
two, then there is always an equivariant isometry
$H^2(S;\ZZ)\cong \Lambda_E\perp \Lambda_E\perp U$ which maps 
to $h$ to $h^+$. In that case the preimages of $e^+$ and $f^+$ are
classes of effective divisors and $\bar S$ is general precisely when
these are the classes of elliptic fibrations. The only way this can fail is
when $S$ contains a $(-2)$-class whose intersection product with
these two isotropic classes is $1$ and $-1$ 
respectively. So this is a property that is detected by the Picard lattice.
\end{none}

We are now ready to set up a period mapping on the level of line
bundles. Let $\LL_E^\times$ be a connected 
component in $\Lambda^-\otimes\CC$ defined by the conditions
$z\cdot z=0$, $z\cdot \bar z>0$ and let $\DD_E$ be its projectivization.
Denote by $\G_E$ the group of isometries of $\Lambda$ that centralize 
$\iota$ and preserve both $h^+$ and $\DD_E$.  
Put $X_E:=\G_E\bs\DD_E$ and let $\Lcal_E$ be the line bundle over
$X_E$ defined by $\LL_E$. Then we have a refined period mapping
\[
G_E\bs (U_E, \eta_E|U_E)\to (X_E,\Lcal_E).
\]
A Torelli theorem asserts that this is an open embedding with image
the complement $X_E^\circ$ of a hypersurface $D_E$ in $X_E$ defined by an
arithmetic arrangement. This arrangement is the collection $\Hcal_E$ of 
hyperplanes
in  $\Lambda^-\otimes\CC$ that are perpendicular to a $(-2)$-vector $r$ in
$\Lambda$ with $r\cdot e^+=1$ and $r\cdot f^+=-1$. This missing
hypersurface parametrizes the special Enriques surfaces. As a subset of 
the Grassmannian, $\Hcal_E$ is an orbit of $\G_E$ and so $D_E$ is
irreducible. The maximal codimension of an intersection of hyperplanes 
from $\Hcal_E$
which meets $\DD_E$ is $2$. Since $Y_E-U_E$ is of codimension $>1$ in 
$Y-E$, it follows from Theorem \ref{thm:gitcon}:

\begin{theorem}\label{thm:enriques} 
The refined period mapping determines an isomorphism of normal projective
 varieties $G_E\bss Y_E^\ss\cong\widehat{X_E^\circ}$ and an 
isomorphism of graded $\CC$-algebra's 
\[
\oplus_{k=0}^\infty H^0(Y_E,\eta^{\otimes k}_E)^{G_E}\cong 
\oplus_{k=0}^\infty H^0(\DD^\circ_E,\Ocal(\LL_E^k))^{\G_E}.
\]
\end{theorem}

A stratification of $G_E\bss Y_E^\ss$ has been worked out by Shah
\cite{shah:enriques}. Sterk shows that if this stratification is transfered
to the arithmetic side of the equation, then it is the one 
that comes with the definition of $\widehat{X_E^\circ}$. He also
proves that the compactication $X_E^{\Sigma (\Hcal_E)}$ of $X_E$ defined by
$\Hcal_E$  is obtained as GIT quotient of a natural blowup of $Y_E$.

\begin{remark}
The singular divisors parametrized by $Y_E$ make up a hypersurface.
This hypersurface is given by a section of a power of $\eta_E$.
On the arithmetic side it must be given by a meromorphic automorphic 
form whose divisor is supported by the arithmetic arrangement that is the
union of $\Hcal_E$ and the hyperplanes perpendicular to $(-2)$-vectors
in $\Lambda^-$. (This form is not  the one exhibited by 
R.~Borcherds in \cite{borcherds:enriques}, since  we also excluded 
the locus of special Enriques surfaces.) We see from this discussion 
that this form already lives on the $18$-dimensional domain of digonal 
$K3$-surfaces of degree $4$.
\end{remark}

\section{Application to smoothing components of triangle singularities}
In this section we give an application in the same spirit
to the semi-universal deformation of a triangle singularity.
Let us first recall what these are. Given integers 
$2\le p_1\le p_2\le p_3$ with $p_1^{-1}+p_2^{-1}+p_3^{-1}<1$, 
then a triangle singularity of type $(p_1,p_2,p_3)$
is a normal surface singularity which admits a (nonminimal) 
resolution consisting of a `central' smooth rational curve
of selfintersection $-1$ and three smooth rational curves
of selfintersection $-p_1,-p_2,-p_3$, which are pairwise disjoint 
and meet the central curve transversally in a single point.
These conditions characterize the singularity up to isomorphism.
This singularity has a $\CC^\times$-action with positive weights:
that is, it sits naturally on a normal affine surface with  $\CC^\times$-action
such that all $\CC^\times$-orbits have the singularity in their
closure. Its semi-universal deformation inherits this $\CC^\times$-action
but here all but one of the weights are positive so that the 
positive weight part is of codimension one in the full
semi-universal deformation. We can realize this singularity 
and the positive weight part of semi-universal deformation
in a projective setting as follows. 

By a standard procedure we can canonically complete this surface to a normal
projective surface $Z_o$ that is smooth away from the $T_{p_1,p_2.p_3}$-singularity.
What has been added at infinity is a $T_{p_1,p_2,p_3}$-curve $T_o$: that is 
a union of smooth rational curves whose dual intersection graph
is a $T$-shape whose arms have  $p_1-1$, $p_2-1$, $p_3-1$ edges. The naturality
guarantees that the $\CC^\times$ action extends to $Z_o$. 
One can further show that $Z_o$ is regular ($H^1(\Ocal_{Z_o})=0$) and that
its dualizing sheaf is trivial with $\CC^\times$-acting on $H^0(\omega_{Z_o})$ 
with degree one. 
A deformation of $Z_o$ which preserves $T_o$ induces a deformation
of the singularity $(Z_o,x_o)$ of positive weight. To be precise,
there exist
\begin{enumerate}
\item[---]a proper flat morphism $Z\to S$
whose fibers are normal surfaces and whose base is an affine scheme, 
\item[---]a distinguished point $o\in S$ such that
the fiber $Z_o$ is identified with its namesake above,
\item[---]an $S$-embedding of $T_o\times S$ in the open subset of $Z$
where $\pi$ is smooth and which over $o$ is the inclusion,  
\item[---]a good $\CC^\times$-action on $\pi$ preserving these data, such
that the local-analytic germ of $\pi$ at $x_o$ is the  
the positive weight part of a semi-universal deformation of
the singularity $(Z_o,x_o)$.
\end{enumerate}
In this situation we can find a generating section $\zeta$ of the
relative dualizing sheaf $\omega_{Z/S}$ such that $t\in\CC^\times$
sends $\zeta$ to $t\zeta$. So for every fiber $s\in S$,
we have a triple $(Z_s,T_o\hookrightarrow  Z_s,\zeta_s)$, 
where $Z_s$ is a regular
normal surface of zero irregularity, $T_o\hookrightarrow Z_s$ an 
embedding in the smooth part of $Z_s$ and
$\zeta_s$ a generator of its dualizing sheaf.  A fiber with only
rational double points as singularities is necessarily a $K3$-surface.
We denote the corresponding subset of $S$ by $S^\circ$. It is open in
$S$. An irreducible component of $S$ which meets $S^\circ$ contains
points with a smooth fiber and is therefore called a \emph{smoothing
component} of $Z_o$. The fibers parametrized by $S-S^\circ$ have a
unique minimally elliptic singularity, other singularities being rational double points.
Denote by $S^\circ\subset S$ the locus over which the only singularities
are rational double points. It is known that $S^\circ$ has pure dimension 
$22-\sum_i p_i$. We shall describe the normalization of the closure of $S^\circ$
in terms of our construction.

\begin{none}
Let  $Q_{p_1,p_2,p_3}$ denote the abstract lattice spanned by the
irreducible components of $T_o$ and with symmetric bilinear form dictated by
the intersection matrix of these irreducible components. This is a nondegenerate
hyperbolic lattice of signature $(1, \sum_i (p_i-1))$.
Any $K3$-surface $S$ that we get as a fiber comes with an embedding 
$T_o\hookrightarrow S$. This clearly induces an injection of  lattices 
$Q_{p_1,p_2,p_3}\hookrightarrow H^2(S;\ZZ)$ whose image lands 
in $\pic (S)$. From this we get that we must have $p_1+p_2+p_3\le 22$.

But the existence of an isometric embedding of $Q_{p_1,p_2,p_3}$ in 
the Picard group of a $K3$-surface need not imply
the existence of an embedding of $T_o$ in that surface.
Let $k(p_1,p_2,p_3)$ be  $6,7,8$ according to whether 
$(p_1,p_2)=(3,3)$ (and hence $p_3\ge 4$)  resp.\ 
$(p_1,p_2)=(2,4)$ (and hence $p_3\ge 5$)  resp.\ 
$(p_1,p_2)=(2,3)$ (and hence $p_3\ge 7$).
Then $Q_{p_1,p_2,p_3}$ contains in an evident manner the affine root 
lattice $\widehat E_{k(p_1,p_2,p_3)}(-1)$. This lattice is negative
semidefinite and  its radical $I(p_1,p_2,p_3)$ is generated by a positive 
integral linear combination
of its basis vectors, called in \cite{looij:triangle2} the
\emph{fundamental  isotropic element} of $Q_{p_1,p_2,p_3}$. Let us denote 
it by $n$.  Suppose that $S$ is a $K3$-surface and that we are given an embedding
of lattices $j: Q_{p_1,p_2,p_3}\hookrightarrow \pic (S)$. 
The theory of  $K3$-surfaces tells us that by composing $j$ with
an isometry of $H^2(S,\ZZ)$ which is plus or minus the identity on the
orthogonal complement of $\pic (Y)$ in $H^2(Y,\ZZ)$, we can arrange that
the fundamental isotropic element is mapped to the class of
an elliptic fibration $Y\to \PP^1$ with each basis vector of
$Q_{p_1,p_2,p_3}$  mapping to the class of a (unique) effective divisor. 
Each of these effective divisors is supported by a configuration
of $(-2)$-curves of type $A$, $D$ or $E$. The 
question is whether we can arrange this in such a manner that these
divisors are all irreducible, for if that is the case, then we clearly
have an embedding  of $T_o$ in $Y$. It is certainly true that in the given 
situation the effective divisors associated to the basis vector of  
$\widehat E_{k(p_1,p_2,p_3)}(-1)$  are supported 
by a single fiber of the elliptic fibration.
One of the main results of \cite{looij:triangle2} states that what we
want is possible precisely when the basis vectors of $\widehat E_k$ are mapped 
to classes of irreducible curves whose union makes up
an entire fiber of the elliptic fibration $Y\to \PP^1$. 
If this fails, then by the Kodaira classification of such fibers, 
the fiber in question must have type $\widehat E_l(-1)$ for some 
$l>k(p_1,p_2,p_3)$. So this can only happen if $k(p_1,p_2,p_3)\in\{
6,7\}$. 

This property can be expressed in purely lattice theoretic terms.  
Let us abbreviate $k$ for $k(p_1,p_2,p_3)$ and $I$ for $I(p_1,p_2,p_3)$.
We call an embedding $j$ of $Q_{p_1,p_2,p_3}$ in an arbitrary
lattice $M$ a \emph{critical} if there exists a $e\in M$ 
with $e\cdot e=-2$, $e\perp j(I)$ and such that $e$ is
neither contained in $j(\widehat E_k)$ nor in its orthogonal
complement; we call $j$ \emph{good} if no such $e$ exists in
in the primitive hull of the image of $j$.
 If $M$ happens to be nondegenerate of hyperbolic signature
$(1,\rho)$ (such as is the case for $\pic (Y)$ above), then $I^\perp/I$ is
a negative definite lattice. The set of $(-2)$-vectors herein
decomposes naturally into irreducible root systems of type $A$, $D$ or
$E$, and one of these contains the image of $j(\widehat E_k)$. So it must be
of type $E_{k(j)}$ for some $k(j)\ge k$. To say that
the embedding is critical is to say that $k(j)>k$. It is clear that this
never happens if $k=8$.
\end{none}

\begin{none}
Still following \cite{looij:triangle2}, this leads to the
following setup.  Let $\LL_{p_1,p_2,p_3}^\times $ be the set
of $(j,z)\in \Hom (Q_{p_1,p_2,p_3},\Lambda )\times(\Lambda\otimes \CC)$
with $j$ a good embedding and $z$ perpendicular to the image of $j$ and
satisfying  the familiar period  conditions $z\cdot z=0$ and $z\cdot\bar z>0$.
The corresponding subset $\DD _{p_1,p_2,p_3}$ of  
$\Hom (Q_{p_1,p_2,p_3},\Lambda )\times\PP(\Lambda\otimes \CC)$
is a disjoint union of type IV domains of dimension $22-\sum_ip_i$.
The orthogonal group $\Orth (\Lambda)$, which operates in an obvious manner
on this space, has  has only finitely many orbits in the connected component set
$\pi_0(\LL_{p_1,p_2,p_3}^\times )$. 
Moreover, the stabilizer of every connected component is 
an arithmetic group acting as such on that component.

Let us say that a pair $(j,[z])\in \DD _{p_1,p_2,p_3}$ is critical if
the embedding defined by $j$ in the orthogonal complement of $[z]$ in 
$\Lambda$ is. Denote by $\DD _{p_1,p_2,p_3}^\circ$ the subset of noncritical
pairs. The following is proved in \cite{looij:triangle2}:
\begin{enumerate}
\item[---]The critical pairs define an arithmetic arrangement in (each
connected component of)  $\DD _{p_1,p_2,p_3}$.
\item[---]The arrangement in question is empty when $k=8$. The maximal
codimension of an intersection from this arrangement is $\le 1$ when
$k=7$ and $\le 2$ when $k=6$. 
\item[---]The quotient $\Orth (\Lambda )\bs (\LL_{p_1,p_2,p_3}^\times)^\circ$ 
is naturally identified with the set $S^\circ$ of $s\in
S$ for which each singular point of $Z_s$ is a rational double point.
\item[---]The dimension of the boundary of $S^\circ$ in $S$  is $\le 10-k$.
\end{enumerate}
\end{none}

\begin{corollary} 
Suppose that $p_1+p_2+p_3-k\le 11$. Then the algebra of meromorphic 
automorphic forms
\[
\oplus_{d\in\ZZ} H^0(\DD_{p_1,p_2,p_3}^\circ,
\Ocal (\LL_{p_1,p_2,p_3}^d))^{\Orth (\Lambda)}
\]
is finitely generated with positive degree generators. Its spectrum is 
the normalization of  $S^\circ$ in $S$.
\end{corollary}
\begin{proof}
Denote by $S_1$ the relevant smoothing component of $S$. 
The elements of degree $d$ of this algebra can be regarded
as functions on $S^\circ$ that are homogeneous of
degree $d$ and so the spectrum in question is the affine completion of
$S^\circ$. The assumption $p_1+p_2+p_3-k\le 11$ implies (according to
the above) that the boundary of $S^\circ$ in $S$ is of codimension $\ge 2$ in 
the closure of $s^\circ$.
So the affine completion of $S^\circ$ is the normalization of $S^\circ$ in $S$.
\end{proof}

A triangle singularity of type  $(p_1,p_2,p_3)$ has a semi-universal
deformation with (irreducible) smooth base if $p_1+p_2+p_3$ is not too large (so
that $\dim \DD_{p_1,p_2,p_3}$ is not too small). For instance, 
$p_1+p_2+p_3\le 12+k$ will do, for in these cases we get a hypersurface.
So then  not only all good embeddings of
$Q_{p_1,p_2,p_3}$ in the $K3$-lattice are conjugate,  but more striking 
perhaps is the conclusion that the corresponding  algebra of meromorphic
automorphic forms 
is a polynomial algebra, a property which seems very hard to establish 
without such a geometric translation. In these cases, the discriminant hypersurface
in the base must be defined by a meromorphic automorphic form
whose divisor is supported by the union of the arrangement above and
the arrangement of hyperplanes in perpendicular to $(-2)$-vectors
that contain the image of $j$. It should have an infinite product expansion.

\end{document}